# LOCAL QUASI-LIKELIHOOD WITH A PARAMETRIC GUIDE[1]

By Jianqing Fan, Yichao Wu and Yang Feng

*Princeton University, North Carolina State University and Princeton University*

Generalized linear models and the quasi-likelihood method extend the ordinary regression models to accommodate more general conditional distributions of the response. Nonparametric methods need no explicit parametric specification, and the resulting model is completely determined by the data themselves. However, nonparametric estimation schemes generally have a slower convergence rate such as the local polynomial smoothing estimation of nonparametric generalized linear models studied in Fan, Heckman and Wand [*J. Amer. Statist. Assoc.* **90** (1995) 141–150]. In this work, we propose a unified family of parametrically-guided nonparametric estimation schemes. This combines the merits of both parametric and nonparametric approaches and enables us to incorporate prior knowledge. Asymptotic results and numerical simulations demonstrate the improvement of our new estimation schemes over the original nonparametric counterpart.

**1. Introduction.** As an extension of the ordinary linear model, the generalized linear model (GLM) broadens techniques of ordinary linear regression to accommodate more general conditional distributions of the response. It was first introduced by Nelder and Wedderburn (1972). Its estimation is based on the iteratively reweighed least squares (IRLS) algorithm, which only requires a relationship between conditional mean and variance instead of its full conditional distribution. This feature was noticed by Wedderburn (1974). In this important further extension, Wedderburn replaced the log-likelihood by a quasi-loglikelihood function. This is usually referred to as the quasi-likelihood method (QLM).

Received October 2008; revised April 2009.

[1]This research was supported in part by NIH Grant R01-GM07261, NSF Grant DMS-07-04337 and DMS-09-05561.

In generalized linear models (GLMs) [McCullagh and Nelder (1989)], a typical parametric assumption is that a transformation of the conditional mean, referred to as the link function, belongs to some parametric family (say, linear or quadratic in the predictor variables). However, misspecification of the parametric family can lead to a completely wrong picture of the underlying conditional mean function. This deficiency of parametric modeling has long been realized in ordinary regression and applies to GLMs as well. It calls for an extension of nonparametric regression techniques to the GLMs. Green and Yandell (1985), O'Sullivan, Yandell and Raynor (1986), and Cox and O'Sullivan (1990) studied the extension to smoothing splines. Tibshirani and Hastie (1987) based their generalization on the "running lines" smoother. Fan, Heckman and Wand (1995) extended the local polynomial fitting technique and includes Staniswalis (1989) as a special case.

Local polynomial smoothing is a useful technique to explore unknown structure in regression and dates back to Stone (1975, 1977). This area blossomed when Fan (1993) provided a deep theoretic understanding and discovered its elegant properties including the automatic boundary correction. Here we focus on local polynomial techniques although the idea can be extended to other nonparametric methods.

Nonparametric methods need no explicit specification of the form of the conditional mean for ordinary regression, and more generally, the link transformation of the conditional mean in the context of GLMs. However, they have in general a slower rate of convergence. In practice, prior knowledge or exploratory studies may provide us some prior information about the shape of the link transformation of the conditional mean. This information is ready to guide us in the nonparametric modeling process. In the literature, parametrically-guided nonparametric estimation methods were proposed to improve over its nonparametric counterpart in the context of density estimation [Hjort and Glad (1995); Naito (2004)] and least squares regression [Glad (1998); Martins-Filho, Mishra and Ullah (2008)]. The idea is very easy to explain in the least squares regression case. Assume that the response $Y$, given a covariate $X$, has a conditional mean $m(x) = E(Y|X = x)$. Once a parametric estimator $m(x, \hat{\beta})$ of $m(x)$ is obtained, any nonparametric method can be applied on $\{Y_i/m(X_i, \hat{\beta}), i = 1, 2, \ldots, n\}$ and $\{Y_i - m(X_i, \hat{\beta}), i = 1, 2, \ldots, n\}$ to estimate $m(x)/m(x, \hat{\beta})$ and $m(x) - m(x, \hat{\beta})$, respectively. The corresponding two final estimators are given by the product of $m(x, \hat{\beta})$ and the nonparametric estimator of $m(x)/m(x, \hat{\beta})$, which serves as a nonparametric correction of the parametric estimator $m(x, \hat{\beta})$, and the sum of $m(x, \hat{\beta})$ and the nonparametric estimator of $m(x) - m(x, \hat{\beta})$, respectively. Theoretically these two parametrically-guided estimators are shown to achieve bias reduction compared to the original nonparametric estimator when $m(\cdot)$ can be approximated by the family $\{m(\cdot, \beta)\}$.

Due to its nice property of bias reduction, it is desirable to extend this parametrically-guided estimation scheme to GLMs and QLM. However, for response with a general distribution other than normal, the regressands $Y/m(X, \hat{\beta})$ and $Y - m(Y, \hat{\beta})$ do not have a nice statistical property to facilitate estimating $m(x)/m(x, \hat{\beta})$ and $m(x) - m(x, \hat{\beta})$ to make the straightforward extension possible. In this work, we take on this problem and propose a unified family of parametrically-guided estimation schemes for QLM. Asymptotic theory and numerical simulations are used to justify our proposed methods. In the literature, similar approaches have been used to reduce variance. Cheng, Peng and Wu (2007) proposed to form a linear combination of a preliminary estimator to reduce variance in smoothing, and Cheng and Hall (2003) studied variance reduction in nonparametric surface estimation.

The rest of the paper is organized as follows. Section 2 presents a fundamental framework of GLMs and QLM. A unified family of parametrically-guided nonparametric estimation schemes is introduced in Section 3. Asymptotic properties are developed to show their improvement over the original nonparametric counterpart in Section 4. Section 5 discusses how to select one parameter in the unified family. Section 6 gives a general pre-asymptotic bandwidth selector based on bias-variance tradeoff. Simulations in Section 7 and real data analysis in Section 8 show our new schemes' finite sample performance in comparison to the original nonparametric method. We conclude with a short discussion in Section 9. Technical proofs are given in the Appendix.

## 2. GLMs and quasi-likelihood models.

Let $(\mathbf{X}_1, Y_1), \ldots, (\mathbf{X}_n, Y_n)$ be a set of i.i.d. random pairs where for each $i$, $Y_i$ is a scalar response variable, and $\mathbf{X}_i$ denotes its corresponding $d$-dimension explanatory covariates having density $f_{\mathbf{X}}$ with support $\mathrm{supp}(f_{\mathbf{X}}) \subseteq \mathbb{R}^d$. In GLMs, we assume that the response's conditional distribution belongs to a one-parameter exponential family

$$(2.1) \qquad f_{Y|\mathbf{X}}(y|\mathbf{x}) = \exp([y\theta(\mathbf{x}) - b(\theta(\mathbf{x}))]/a(\phi) + c(y, \phi)),$$

where $a(\cdot)$, $b(\cdot)$ and $c(\cdot, \cdot)$ are some known functions, $\phi$ is the dispersion parameter and $\theta$ is the canonical parameter. For (2.1), the response has conditional mean $\mu(\mathbf{x}) = b'(\theta(\mathbf{x}))$ and conditional variance $\mathrm{var}(Y|\mathbf{X} = \mathbf{x}) = a(\phi)b''(\theta(\mathbf{x}))$.

Parametric GLMs assume that $\eta(\mathbf{x}) = g(\mu(\mathbf{x}))$ and $\eta(\mathbf{x}) = \beta_0 + \mathbf{x}^T \boldsymbol{\beta}$ for some monotonic link $g(\cdot)$. When the canonical link $g = (b')^{-1}$ is used, the composition $g \circ b'(\cdot)$ reduces to the identity function and $\theta(\mathbf{x}) \equiv \eta(\mathbf{x})$. In this case, (2.1) simplifies to $f_{Y|\mathbf{X}}(y|\mathbf{x}) = \exp([y\eta(\mathbf{x}) - b(\eta(\mathbf{x}))]/a(\phi) + c(y, \phi))$.

In common practice, the full likelihood may be unavailable. However, the relationship between the conditional mean and variance may be readily available. In this case, estimation of $\mu(\mathbf{x})$ can be achieved by replacing the conditional log-likelihood $\log f_{Y|\mathbf{X}}(y|\mathbf{x})$ by a quasi-log-likelihood function $Q(\mu(\mathbf{x}), y)$. When we assume that $\mathrm{var}(Y|\mathbf{X} = \mathbf{x}) = V(\mu(\mathbf{x}))$ for some known positive variance function $V(\cdot)$, the corresponding $Q(\mu, y)$ satisfies

$$(2.2) \qquad U(w, y) = \frac{\partial}{\partial w} Q(w, y) = \frac{y - w}{V(w)}.$$

More explicitly, $Q(\mu, y) = \int_y^\mu (y - w)/V(w)\,dw$. For more details on QLM, see Wedderburn (1974) and Chapter 9 of McCullagh and Nelder (1989). The quasi-score (2.2) possesses properties similar to those of the usual log-likelihood score function. Note that the loglikelihood of (2.1) is a special case of quasi-likelihood function with $V(\cdot) = a(\phi)b'' \circ (b')^{-1}(\cdot)$.

Due its generality, we will focus on QLM. Fan, Heckman and Wand (1995) introduced nonparametric QLM by extending the local polynomial techniques. We will follow their framework and notation. To ease our presentation, we focus on the one-dimension case as the extension to the multivariate case is straightforward. For the one-dimension case, our data consist of $n$ pairs of observations $\{(X_i, Y_i), i = 1, 2, \ldots, n\}$.

To enhance flexibility, Fan, Heckman and Wand (1995) modeled $\eta(x)$ nonparametrically. For any $x_0$ in its domain, the local polynomial estimator of $\eta(x_0)$ is given by $\hat{\eta}(x_0) \equiv \hat{\eta}(x_0; p, h) = \hat{\beta}_0$ where $\hat{\boldsymbol{\beta}} = (\hat{\beta}_0, \hat{\beta}_1, \ldots, \hat{\beta}_p)^T$ maximizes the locally weighted quasi-likelihood function

$$(2.3) \qquad Q(\boldsymbol{\beta}) = Q(\boldsymbol{\beta}; h, x_0) = \sum_{i=1}^n Q(g^{-1}(\mathbf{X}_i^T \boldsymbol{\beta}), Y_i) K_h(X_i - x_0),$$

where, with slight abuse of notation, we define $\mathbf{X}_i = (1, X_i - x_0, \ldots, (X_i - x_0)^p)^T$ and $\boldsymbol{\beta} = (\beta_0, \beta_1, \ldots, \beta_p)^T$. Whenever there is no confusion, the extra arguments are dropped and $Q(\boldsymbol{\beta})$ is used, similarly for some other notation. Here $p$ is the order of local polynomial fitting and $K_h(\cdot) = K(\cdot/h)/h$ is a re-scaling of the kernel function $K(\cdot)$ with a smoothing bandwidth $h$.

## 3. Nonparametric quasi-maximum likelihood with a parametric guide.
As argued in the introduction, prior knowledge, physical model or exploratory analysis may give us some useful information that $\eta(x)$ falls approximately into a parametric family $\{\eta(x, \boldsymbol{\alpha}) : \boldsymbol{\alpha} = (\alpha_1, \alpha_2, \ldots, \alpha_q)^T \in \mathbb{A} \subset \mathbb{R}^q\}$. In this section, we present a family of estimation schemes by incorporating the available useful shape information of $\eta(x)$ to guide us while estimating $\eta(x)$. Within the parametric family $\eta(x, \boldsymbol{\alpha})$, we find the optimal fit by maximizing

$$(3.1) \qquad \sum_{i=1}^n Q(g^{-1}(\eta(X_i, \boldsymbol{\alpha})), Y_i)$$

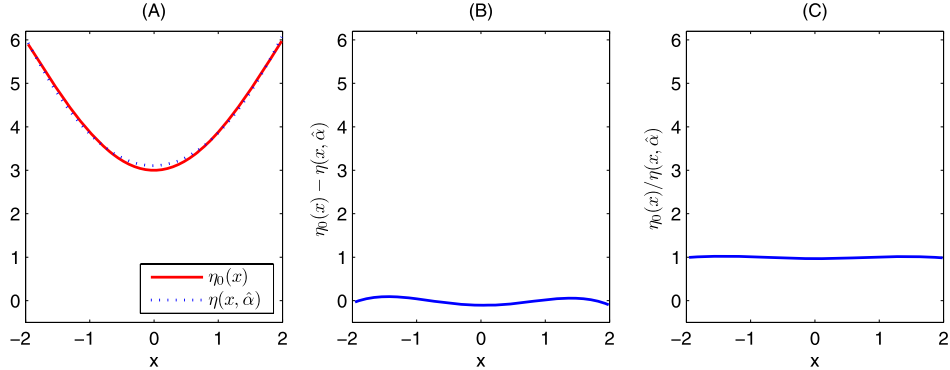

Fig. 1. *Plots of true $\eta(\cdot)$, estimated guide $\eta(\cdot, \hat{\boldsymbol{\alpha}})$, difference $\eta(\cdot) - \eta(\cdot, \hat{\boldsymbol{\alpha}})$ and ratio $\eta(\cdot)/\eta(\cdot, \hat{\boldsymbol{\alpha}})$ for one random sample in Example 7.1.*

with respect to $\boldsymbol{\alpha} \in \mathbb{A}$. Denote the best fit by $\eta(x, \hat{\boldsymbol{\alpha}})$ where $\hat{\boldsymbol{\alpha}}$ is the maximizer of (3.1).

3.1. *Bias reduction.* In the local polynomial fitting framework, the bias is due to the approximation error of the Taylor expansion. The smaller approximation error the less bias in the local polynomial estimator. Recall that we identify some parametric family $\{\eta(x, \boldsymbol{\alpha}): \boldsymbol{\alpha} \in \mathbb{A}\}$ based on exploratory studies or prior knowledge and find the best fit $\eta(x, \hat{\boldsymbol{\alpha}})$ within this family. As a result, $\eta(x, \hat{\boldsymbol{\alpha}})$ should capture the major shape of $\eta(x)$ and consequently $\eta(x)/\eta(x, \hat{\boldsymbol{\alpha}})$ and $\eta(x) - \eta(x, \hat{\boldsymbol{\alpha}})$ have less variation (smoother) than the original $\eta(x)$ does. Consequently, they are easier to be approximated and the approximation errors in their corresponding Taylor expansions are smaller than those of the original function $\eta(x)$. For example, the true $\eta(\cdot)$ is given by $\eta_0(x) = 3\sin(\frac{\pi}{4}x - \frac{\pi}{2}) + 6$ for $x \in [-2, 2]$ [as shown by the solid line in panel (A) of Figure 1] in our Poisson simulation Example 7.1. Nonparametric estimate $\hat{\eta}(\cdot)$ is given by the dotted line in Figure 2 and indicates a parabolic shape. Hence we identify a parametric family, $\{\eta(x, \boldsymbol{\alpha}) = \alpha_1 + \alpha_2 x + \alpha_3 x^2: \boldsymbol{\alpha} = (\alpha_1, \alpha_2, \alpha_3)^T \in \mathbb{R}^3\}$, within which the best fit is given by the dotted line in panel (A) of Figure 1. The difference $\eta(x) - \eta(x, \hat{\boldsymbol{\alpha}})$ and ratio $\eta(x)/\eta(x, \hat{\boldsymbol{\alpha}})$ are shown in panels (B) and (C) of Figure 1, respectively. We can see that the difference and ratio functions are much flatter than the original function $\eta(\cdot)$ as desired.

Based on the above argument, two different estimation schemes corresponding to multiplicative and additive corrections are introduced in Sections 3.2 and 3.3, respectively. They are special cases of a unified family of corrections presented in Section 3.4.

### 3.2. *Multiplicative correction.*

Consider the multiplicative identity

$$\eta(x) \equiv \eta(x, \boldsymbol{\alpha}) r_m(x),$$

where $r_m(x) = \eta(x)/\eta(x, \boldsymbol{\alpha})$. When $\eta(x, \hat{\boldsymbol{\alpha}})$ is a good estimate, the ratio $r_m(x)$ becomes almost flat and allows the choice of a larger bandwidth. For any $x_0$, we may estimate $r_m(x_0)$ by maximizing local quasi-likelihood

$$\sum_{i=1}^{n} Q(g^{-1}(\mathbf{X}_i^T \boldsymbol{\beta} \eta(X_i, \hat{\boldsymbol{\alpha}})), Y_i) K_h(X_i - x_0)$$

with respect to $\boldsymbol{\beta}$ and set $\hat{r}(x_0) = \hat{\beta}_0$, the first component of the maximizer $\boldsymbol{\beta}$. Then $\eta(x_0)$ can be estimated by $\eta(x_0, \hat{\boldsymbol{\alpha}})\hat{r}(x_0)$. This two-step formulation is equivalent to the following one-step estimation.

Locally approximating $r_m(\cdot)$ by a polynomial function and re-scaling it by a factor $\eta(x_0, \hat{\boldsymbol{\alpha}})$, we have the local quasi-likelihood

$$(3.2) \quad \begin{aligned} Q_m(\boldsymbol{\beta}) &\equiv Q_m(\boldsymbol{\beta}; h, x_0, \hat{\boldsymbol{\alpha}}) \\ &= \sum_{i=1}^{n} Q(g^{-1}(\mathbf{X}_i^T \boldsymbol{\beta} \eta(X_i, \hat{\boldsymbol{\alpha}})/\eta(x_0, \hat{\boldsymbol{\alpha}})), Y_i) K_h(X_i - x_0). \end{aligned}$$

We maximize (3.3) with respect to $\boldsymbol{\beta}$, and set the final estimator $\hat{\eta}_m(x_0) \equiv \hat{\eta}_m(x_0; p, h, \hat{\boldsymbol{\alpha}}) = \hat{\beta}_0$. In this formulation, the Taylor expansion $\mathbf{X}_i^T \boldsymbol{\beta}$ is supposed to approximate $\eta(X_i)\eta(x_0, \hat{\boldsymbol{\alpha}})/\eta(X_i, \hat{\boldsymbol{\alpha}})$ locally at $x = x_0$. This immediately justifies setting $\hat{\beta}_0$ as our estimator.

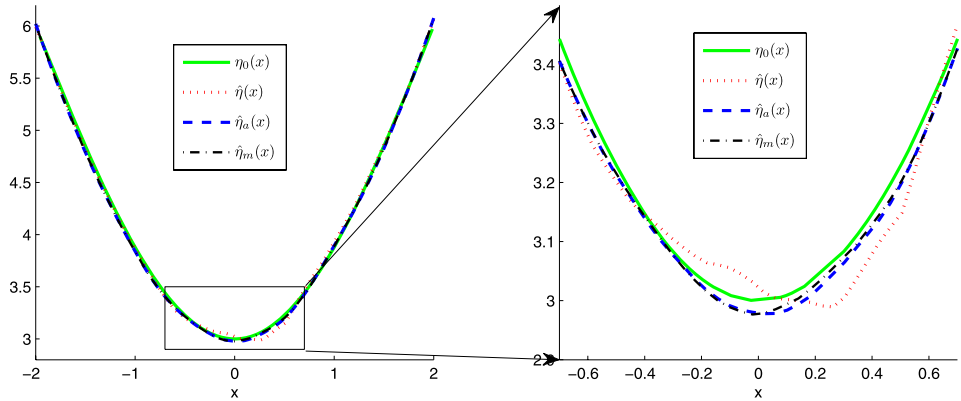

FIG. 2. *Plots of true $\eta(\cdot)$ nonparametric estimate $\hat{\eta}(\cdot)$ and two parametrically-guided estimates $\hat{\eta}_a(\cdot)$ and $\hat{\eta}_m(\cdot)$ for one random sample in Example 7.1 are shown on the left panel. A zoom-in view of the squared region is given on the right panel.*

### 3.3. *Additive correction.* The other additive identity

$$\eta(x) \equiv \eta(x, \boldsymbol{\alpha}) + r_a(x)$$

with $r_a(x) = \eta(x) - \eta(x, \boldsymbol{\alpha})$ leads to another parametrically-guided nonparametric estimator $\hat{\eta}_a(x_0) \equiv \hat{\eta}_a(x_0; p, h, \hat{\boldsymbol{\alpha}}) = \hat{\beta}_0$ where $\hat{\beta}_0$ is the first component of the maximizer of

$$
\begin{aligned}
(3.3) \quad Q_a(\boldsymbol{\beta}) &\equiv Q_a(\boldsymbol{\beta}; h, x_0, \hat{\boldsymbol{\alpha}}) \\
&= \sum_{i=1}^{n} Q(g^{-1}(\eta(X_i, \hat{\boldsymbol{\alpha}}) - \eta(x_0, \hat{\boldsymbol{\alpha}}) + \mathbf{X}_i^T \boldsymbol{\beta}), Y_i) K_h(X_i - x_0)
\end{aligned}
$$

with respect to $\boldsymbol{\beta}$. Similarly, the expansion $\mathbf{X}_i^T \boldsymbol{\beta}$ in this formulation targets at approximating $\eta(X_i) - \eta(X_i, \hat{\boldsymbol{\alpha}}) + \eta(x_0, \hat{\boldsymbol{\alpha}})$ locally at $x = x_0$. Hence $\hat{\beta}_0$ estimates $\eta(x_0)$.

### 3.4. *A unified family of corrections.* As in Martins-Filho, Mishra and Ullah (2008), we consider a more general identity $\eta(x) \equiv \eta(x, \boldsymbol{\alpha}) + r_u(x)\eta(x, \boldsymbol{\alpha})^\gamma$ with $r_u(x) = (\eta(x) - \eta(x, \boldsymbol{\alpha}))/\eta(x, \boldsymbol{\alpha})^\gamma$ for some $\gamma \geq 0$, we can estimate $\eta(x_0)$ by $\hat{\eta}_u(x_0) = \eta(x_0, \hat{\boldsymbol{\alpha}}) + \hat{r}_u(x_0)\eta(x_0, \hat{\boldsymbol{\alpha}})^\gamma$. Here $\hat{r}_u(x_0)$ is given by the first component $\hat{\beta}_0$ of the maximizer of

$$\sum_{i=1}^{n} Q(g^{-1}(\eta(X_i, \hat{\boldsymbol{\alpha}}) + \mathbf{X}_i^T \boldsymbol{\beta}\eta(X_i, \hat{\boldsymbol{\alpha}})^\gamma), Y_i) K_h(X_i - x_0).$$

As in Section 3.2, an equivalent one-step estimation is available. Let $\hat{\beta}_0$ be the first component of the maximizer of

$$
\begin{aligned}
(3.4) \quad Q_u(\boldsymbol{\beta}) &\equiv Q_u(\boldsymbol{\beta}; h, x_0, \hat{\boldsymbol{\alpha}}) \\
&= \sum_{i=1}^{n} Q(g^{-1}(\eta(X_i, \hat{\boldsymbol{\alpha}}) \\
&\qquad\qquad + (\mathbf{X}_i^T \boldsymbol{\beta} - \eta(x_0, \hat{\boldsymbol{\alpha}}))\eta(X_i, \hat{\boldsymbol{\alpha}})^\gamma / \eta(x_0, \hat{\boldsymbol{\alpha}})^\gamma), Y_i) \\
&\qquad \times K_h(X_i - x_0)
\end{aligned}
$$

with respect to $\boldsymbol{\beta}$. Then $\hat{\beta}_0$ directly estimates $\eta(x_0)$ and is the same as $\hat{\eta}_u(x_0)$. We prefer (3.4) since it facilitates our theoretical development.

Note that this unified estimator includes the additive and multiplicative corrections as special cases by setting $\gamma = 0$ and 1, respectively.

## 4. Asymptotic properties.

We assume that our data $\{(X_i, Y_i), i = 1, 2, \ldots, n\}$ are generated from the quasi-likelihood model with unknown true $\eta_0(x)$. Asymptotic properties of our final estimates are achieved in two steps: establish asymptotic properties with a fixed parametric guide in Section 4.1 and show that the same asymptotic properties apply to the case with an estimated parametric guide in Section 4.2.

In the univariate case, we denote the marginal density of $X$ by $f_X$. The asymptotic properties for local polynomial estimator are different for $x_0$ lying in the interior of $\mathrm{supp}(f_X)$ from for $x_0$ lying near the boundary. Suppose that $K$ is supported on $[-1, 1]$. Then the support of $K_h(x_0 - \cdot)$ is $\mathcal{E}_{x_0, h} = \{z : |z - x_0| \leq h\}$. We will call $x_0$ an interior point of $\mathrm{supp}(f_X)$ if $\mathcal{E}_{x_0, h} \subset \mathrm{supp}(f_X)$ and a boundary point otherwise. If $\mathrm{supp}(f_X) = [a, b]$, then $x_0$ is a boundary point if and only if $x_0 = a + \alpha h$ or $x_0 = b - \alpha h$ for some $0 \leq \alpha < 1$. Denote $\mathcal{D}_{x_0, h} = \{z : x_0 - hz \in \mathrm{supp}(f_X)\} \cap [-1, 1]$. For any measurable set $\mathcal{A} \subset \mathcal{R}$, define $\nu_l(\mathcal{A}) = \int_{\mathcal{A}} z^l K(z) \, dz$. Let $\mathbf{N}_p(\mathcal{A})$ be the $(p+1) \times (p+1)$ matrix having $(i, j)$ entry equal to $\nu_{i+j-2}(\mathcal{A})$, and let $\mathbf{M}_{r,p}(z; \mathcal{A})$ be the same as $\mathbf{N}_p(\mathcal{A})$, but with the $(r+1)$th column replaced by $(1, z, \ldots, z^p)^T$. Then for $|\mathbf{N}_p(\mathcal{A})| \neq 0$, define

$$(4.1) \qquad K_{r,p}(z; \mathcal{A}) = r! \{|\mathbf{M}_{r,p}(z; \mathcal{A})| / |\mathbf{N}_p(\mathcal{A})|\} K(z).$$

When $[-1, 1] \subseteq \mathcal{A}$, we will suppress $\mathcal{A}$ and simply write $\nu_l, \mathbf{N}_p, \mathbf{M}_{r,p}$, and $K_{r,p}$. It can be shown that $(-1)^r K_{r,p}(\cdot; \mathcal{A})$ is an order $(r, s)$ kernel as defined by Gasser, Müller and Mammitzsch (1985) where $s = p + 1$ if $p - r$ is odd, and $s = p + 2$ if $p - r$ is even. It is an equivalent kernel induced by the local polynomial fitting [Fan and Gijbels (1995)]. This family of kernels is useful for giving concise expressions for the asymptotic distribution of local polynomial estimator for $x_0$ lying either in the interior of $\mathrm{supp}(f_X)$ or near its boundaries. Denote $\rho(x_0) = \{g'(\mu(x_0))^2 V(\mu(x_0))\}^{-1}$. Note that when the model belongs to a one-parameter exponential family and the canonical link is used then $g'(\mu(x_0)) = 1/\mathrm{var}(Y|X = x_0)$, and $\rho(x_0) = \mathrm{var}(Y|X = x_0)$, if the variance function $V(\cdot)$ is correctly specified. The asymptotic variance of our local polynomial estimator depends on

$$\sigma_{r,s,p}^2(x_0; K, \mathcal{A})$$
$$= \mathrm{var}(Y|X = x_0) g'(\mu(x_0))^2 f_X(x_0)^{-1} \int_{\mathcal{A}} K_{r,p}(z; \mathcal{A}) K_{s,p}(z; \mathcal{A}) \, dz.$$

Since the multiplicative and the additive corrections are both special cases of the unified family of corrections, we only consider the asymptotic properties for the unified family of corrections.

### 4.1. Asymptotic properties with a fixed guide.
Recall that our parametrically-guided nonparametric estimators are achieved by maximizing $Q_u(\boldsymbol{\beta}; h, x_0, \hat{\boldsymbol{\alpha}})$ defined by (3.4). Note that the definition of $Q_u(\boldsymbol{\beta}; h, x_0, \hat{\boldsymbol{\alpha}})$ involves $\hat{\boldsymbol{\alpha}}$ which corresponds to the best fit within the parametric family $\{\eta(x, \boldsymbol{\alpha}), \boldsymbol{\alpha} \in \mathbb{A}\}$ and depends on our data $\{(X_i, Y_i), i = 1, 2, \ldots, n\}$. This dependency consequently makes it intractable to directly study the asymptotical properties of the maximizer of $Q_u(\boldsymbol{\beta}; h, x_0, \hat{\boldsymbol{\alpha}})$. To avoid the complication caused by the use of the estimated $\hat{\boldsymbol{\alpha}}$, we first consider the case with a fixed guide $\eta(x, \boldsymbol{\alpha})$.

When a fixed guide $\eta(\cdot, \boldsymbol{\alpha})$ is used, the asymptotic normality of the corresponding estimator $\hat{\eta}_u(x_0; p, h, \boldsymbol{\alpha})$ are given by Theorem 1.

THEOREM 1. *Let* $p > 0$, $\gamma \geq 0$, *and assume that* $h = h_n \to 0$, $nh^{2p+1} \to \infty$, *and* $nh^{2p+3} < \infty$ *as* $n \to \infty$. *Under conditions* (A1)–(A5) *stated in the Appendix, if* $x_0$ *is a fixed point in the interior of* supp($f_X$) *satisfying* $\eta(x_0, \boldsymbol{\alpha}) \neq 0$, *then we have*

$$(4.2) \qquad \frac{\sqrt{nh}}{\sigma_{0,0,p}(x_0; K)} [\hat{\eta}_u(x_0; p, h, \boldsymbol{\alpha}) - \eta_0(x_0) - Bias] \xrightarrow{\mathcal{D}} N(0, 1),$$

*where the bias term is given by* $Bias_o$ *for odd* $p$ *and* $Bias_e$ *for even* $p$ *defined by*

$$Bias_o = \frac{\eta(x_0, \boldsymbol{\alpha})^\gamma}{(p+1)!} \left( \frac{\eta_0(\cdot) - \eta(\cdot, \boldsymbol{\alpha})}{\eta(\cdot, \boldsymbol{\alpha})^\gamma} \right)^{(p+1)} (x_0) h^{p+1}$$
$$\times \left( \int z^{p+1} K_{0,p}(z) \, dz \right) \{1 + O(h)\}$$

*and*

$$Bias_e = \left\{ \int z^{p+2} K_{0,p}(z) \, dz \frac{1}{(p+2)!} \left( \frac{\eta_0(\cdot) - \eta(\cdot, \boldsymbol{\alpha})}{\eta(\cdot, \boldsymbol{\alpha})^\gamma} \right)^{(p+2)} (x_0) \eta(x_0, \boldsymbol{\alpha})^\gamma \right.$$
$$+ \frac{1}{(p+1)!} \left( \frac{\eta_0(\cdot) - \eta(\cdot, \boldsymbol{\alpha})}{\eta(\cdot, \boldsymbol{\alpha})^\gamma} \right)^{(p+1)}$$
$$\left. \times (x_0) \frac{(\rho \eta(\cdot, \boldsymbol{\alpha})^{2\gamma} f_X)'(x_0)}{(\rho \eta(\cdot, \boldsymbol{\alpha})^\gamma f_X)(x_0)} \int z^{p+2} K_{0,p}(z) \, dz \right\}$$
$$\times h^{p+2} \{1 + O(h)\}.$$

*If* $x_0 = x_n$ *is of the form* $x_0 = x_\delta + ch$ *satisfying* $\eta(x_0, \boldsymbol{\alpha}) \neq 0$ *where* $x_\delta$ *is a point on the boundary of* supp($f_X$) *and* $c \in [-1, 1]$, *then* (4.2) *holds with* $\sigma_{0,0,p}^2(x_0; K)$, *and* $\int z^{p+1} K_{0,p}(z) \, dz$ *replaced by* $\sigma_{0,0,p}^2(x_0; K, \mathcal{D}_{x_0,h})$, *and* $\int_{\mathcal{D}_{x_0,h}} z^{p+1} K_{0,p}(z; \mathcal{D}_{x_0,h}) \, dz$.

REMARK 1. Note that we use $\eta(x_0, \boldsymbol{\alpha})$ in the denominator, which poses difficulty handling any zero point of $\eta(\cdot, \boldsymbol{\alpha})$, that is, $x_0$ satisfying $\eta(x_0, \boldsymbol{\alpha}) = 0$. These zero points are ruled out in Theorem 1. Similar observation was made by Hjort and Glad (1995). However, this difficulty does not occur in our limited numerical experiments.

REMARK 2. To simplify our presentation, we only state the asymptotic normality of the estimator for the function $\eta_0(\cdot)$. However, our method can also estimate its high order derivatives, for which the asymptotic properties can be found in Proposition 1 in the Appendix.

### 4.2. *Asymptotic properties with an estimated guide.*

Note that our proposed estimation schemes use the best parametric fit $\eta(x, \hat{\boldsymbol{\alpha}})$ estimated based on our data instead of a fixed guide $\eta(x, \boldsymbol{\alpha})$. Compared to the simpler case with a fixed guide, the variability of parameter estimation now influences the asymptotic result. However, we shall show below that asymptotically there is no precision loss caused by the additional estimation step.

Clearly, the parametric family used in the first step of our estimation schemes is most likely an incorrect specification. Consequently, the first-stage parametric estimator is a maximum quasi-likelihood estimator with a misspecified model. As in Hurvich and Tsai (1995), we denote the proposed parametric joint density of $(X_i, Y_i)$ and the corresponding actual unknown joint density by $f(x, y; \boldsymbol{\alpha}) = f_X(x) \exp(Q(g^{-1}(\eta(x, \boldsymbol{\alpha})), y))$ and $f_0(x, y) = f_X(x) \exp(Q(g^{-1}(\eta_0(x)), y))$, respectively, where $f_X(\cdot)$ is the marginal density of $X$. Denote by $\boldsymbol{\alpha}_0$, the pseudo parameter value that minimizes the Kullback–Liebler distance between $f(x, y; \boldsymbol{\alpha})$ and $f_0(x, y)$, that is,

$$\boldsymbol{\alpha}_0 \stackrel{\triangle}{=} \underset{\boldsymbol{\alpha} \in \mathbb{A}}{\operatorname{argmin}}\, E \log\left( \frac{f_0(X, Y)}{f(X, Y; \boldsymbol{\alpha})} \right)$$

$$\equiv \underset{\boldsymbol{\alpha} \in \mathbb{A}}{\operatorname{argmin}} \int\int (Q(g^{-1}(\eta_0(x)), y) - Q(g^{-1}(\eta(x, \boldsymbol{\alpha})), y)) f_0(x, y)\, dy\, dx,$$

where the expectation $E$ is taken with respect to the unknown true density.

To proceed, we make regularity assumptions (B1)–(B5) given in the Appendix to assure that the pseudo-maximum quasi-likelihood estimator $\hat{\boldsymbol{\alpha}}$ is $\sqrt{n}$-consistent of $\boldsymbol{\alpha}_0$, that is, $\sqrt{n}(\hat{\boldsymbol{\alpha}} - \boldsymbol{\alpha}_0) = O_p(1)$ [see White (1982)].

THEOREM 2. *Under additional conditions* (B1)–(B5), *the asymptotic results, with* $\boldsymbol{\alpha}$ *replaced by* $\boldsymbol{\alpha}_0$, *of Theorem 1 continue to hold when an estimated fit* $\eta(x, \hat{\boldsymbol{\alpha}})$ *is used.*

REMARK 3. Note that our theoretical results include those of the original nonparametric method in Fan, Heckman and Wand (1995) as a special case by setting a constant guide, say $\eta(\cdot, \boldsymbol{\alpha}) = 1$. For our parametrically-guided estimation, the asymptotic bias is determined by both the guide $\eta(\cdot, \boldsymbol{\alpha})$ and $\gamma$. When the same smoothing bandwidth is used, our theoretical results allow straightforward comparison between the original nonparametric method and our new parametrically-guided estimation schemes with different $\gamma$. Although the asymptotic variance remains the same, the advantage of using a parametric guide is that it can reduce the asymptotic bias. For example, when $p = 1$ and $h$ is the same, our parametrically-guided method asymptotically reduces integrated squared bias provided that

$$
(4.3) \quad
\begin{aligned}
\inf_{\gamma \geq 0} \int_{\mathrm{supp}(f_X)} & \left( \eta(x, \boldsymbol{\alpha}_0)^\gamma \left( \frac{\eta_0(\cdot) - \eta(\cdot, \boldsymbol{\alpha}_0)}{\eta(\cdot, \boldsymbol{\alpha}_0)^\gamma} \right)^{(2)}(x) \right)^2 dx \\
& < \int_{\mathrm{supp}(f_X)} (\eta_0^{(2)}(x))^2 \, dx.
\end{aligned}
$$

However, this is only one part of the whole story because when the guide is appropriately selected, our parametrically-guided estimation schemes will select larger smoothing bandwidths (as the correction function is smoother) and, consequently, improve performance by reducing variance as well.

**5. Selection of $\boldsymbol{\gamma}$.** Equation (4.3) can be used as a general rule of thumb for identifying an appropriate parametric guide and selecting $\gamma$ by minimizing

$$
\theta_\gamma = \int_{\mathrm{supp}(f_X)} \left( \eta(x, \boldsymbol{\alpha}_0)^\gamma \left( \frac{\eta_0(\cdot) - \eta(\cdot, \boldsymbol{\alpha}_0)}{\eta(\cdot, \boldsymbol{\alpha}_0)^\gamma} \right)^{(2)}(x) \right)^2 dx,
$$

namely, the quantity on the left-hand side.

In finite-sample applications, we obtain the best fit $\eta(x, \hat{\boldsymbol{\alpha}})$ for each potential parametric guide family $\eta(x, \boldsymbol{\alpha})$ and use local polynomial smoothing to estimate the second order derivative function $(\frac{\hat{\eta}(\cdot) - \eta(\cdot, \hat{\boldsymbol{\alpha}})}{\eta(\cdot, \hat{\boldsymbol{\alpha}})^\gamma})^{(2)}(x)$. Then we define

$$
(5.1) \quad \hat{\theta}_\gamma = \int_{\mathrm{supp}(f_X)} \left( \eta(x, \hat{\boldsymbol{\alpha}})^\gamma \left( \frac{\hat{\eta}(\cdot) - \eta(\cdot, \hat{\boldsymbol{\alpha}})}{\eta(\cdot, \hat{\boldsymbol{\alpha}})^\gamma} \right)^{(2)}(x) \right)^2 dx.
$$

Treating $\hat{\theta}_\gamma$ as a function of the parametric family and $\gamma$, we can find the best parametric guide and its corresponding best $\gamma$ by minimizing $\hat{\theta}_\gamma$.

For the case of least squares regression, Huang and Fan (1999) studied convergence rate of nonparametric estimators of quadratic regression functionals such as the quantities on both sides of (4.3). We can apply their Theorems 4.1–4.4 to get the convergence of our plug-in estimator $\hat{\theta}_\gamma$ by noting that $\hat{\boldsymbol{\alpha}}$ converges to $\boldsymbol{\alpha}_0$ with a faster speed. However, the corresponding theory for the more general GLM and quasi-likelihood method is not available. A serious treatment for this kind of problem is very technical. It requires a full paper to address the issues and is beyond our current scope.

In simulation examples of Section 7, for each example we generate 10 additional samples. Based on these 10 samples, we use the Extended Residual Squares Criterion (ERSC) [see equation (5.6) of Fan, Farmen and Gijbels (1998)] to select the smoothing bandwidth for estimating the second order derivative $\left(\frac{\hat{\eta}(\cdot) - \eta(\cdot, \hat{\boldsymbol{\alpha}})}{\eta(\cdot, \hat{\boldsymbol{\alpha}})^{\gamma}}\right)^{(2)}(x)$ for each $\gamma$. Then we can evaluate $\hat{\theta}_{\gamma}$ for each sample and a grid $\gamma_{j}$, $j = 1, 2, \ldots, J$ of $\gamma$, denoted by $\hat{\theta}_{\gamma_{j}, i}$ for $i = 1, 2, \ldots, 10$ and $j = 1, 2, \ldots, J$. Then we select the minimizer $\hat{j} = \operatorname{argmin}_{j} \sum_{i=1}^{10} \hat{\theta}_{\gamma_{j}, i}$ and set $\gamma_{\hat{j}}$ as the tuned $\gamma$.

From our simulations, the improvement from the additive or multiplicative correction to the best $\gamma$ is much smaller than the improvement from the original method to the additive or multiplicative correction. In other words, the sensitivity of $\gamma$ on the performance improvement is not very high. This is due in part to the choice of the parametric guides which usually capture the main shape. Thus, in application, we suggest that a simple and effective method is to try a few discrete values of $\gamma$ [including additive ($\gamma = 0$) and multiplicative ($\gamma = 1$) guides as specific examples] and to pick the value of $\gamma$ by the cross-validation. This will result in an improved performance over the vanilla nonparametric approach, if that approach is also included in the cross-validation comparison.

**6. Pre-asymptotic bandwidth selection.** While optimizing (2.3) and (3.4), we need to tune the corresponding smoothing bandwidths. In this work, we will use the pre-asymptotic bandwidth selection method introduced in Fan and Gijbels (1995) and Fan, Farmen and Gijbels (1998) which is based on the bias-variance tradeoff.

6.1. *Estimating bias and variance.* Without loss of generality, we use (3.4) to demonstrate the idea. It will include (2.3) as a special case by using a constant guide. In the remainder of this section, we denote $\hat{\boldsymbol{\beta}} = \hat{\boldsymbol{\beta}}(x_0, \hat{\boldsymbol{\alpha}}) = \operatorname{argmax}_{\boldsymbol{\beta}} Q_u(\boldsymbol{\beta}; h, x_0, \hat{\boldsymbol{\alpha}})$. The bias of the estimate $\hat{\boldsymbol{\beta}}$ comes from the approximation error in the Taylor expansion. Denote the approximation error at $X_i$ by

$$r(X_i) = \eta_0(X_i) - \eta(X_i, \hat{\boldsymbol{\alpha}})$$
$$- \frac{\eta(X_i, \hat{\boldsymbol{\alpha}})^{\gamma}}{\eta(x_0, \hat{\boldsymbol{\alpha}})^{\gamma}} \sum_{j=0}^{p} \left(\frac{\eta_0 - \eta(\cdot, \hat{\boldsymbol{\alpha}})}{\eta(\cdot, \hat{\boldsymbol{\alpha}})^{\gamma}}\right)^{(j)}(x_0) \frac{(X_i - x_0)^j}{j!}.$$

Suppose that the $(p + a + 1)$th derivatives of functions $\eta_0(\cdot)$ and $\eta(\cdot, \hat{\boldsymbol{\alpha}})$ exist at $x_0$ for some integer $a > 0$. Further expansions of $\eta_0(X_i)$ and $\eta(X_i, \hat{\boldsymbol{\alpha}})$ give

$$r(X_i) \approx \frac{\eta(X_i, \hat{\boldsymbol{\alpha}})^{\gamma}}{\eta(x_0, \hat{\boldsymbol{\alpha}})^{\gamma}} \sum_{j=1}^{a} \left(\frac{\eta_0 - \eta(\cdot, \hat{\boldsymbol{\alpha}})}{\eta(\cdot, \hat{\boldsymbol{\alpha}})^{\gamma}}\right)^{(p+j)}(x_0) \frac{(X_i - x_0)^{p+j}}{(p+j)!} \triangleq r_i.$$

Here the choice of $a$, the approximation order, will affect the performance of the estimated bias. Practically, it can be chosen as $a = 1$ or $2$.

Now pretend that the approximated approximation errors $r_i$ are known. A more accurate local quasi-log-likelihood is

$$Q_u^*(\boldsymbol{\beta}) = Q_u^*(\boldsymbol{\beta}; h, x_0, \hat{\boldsymbol{\alpha}}) = \sum_{i=1}^n Q(g^{-1}(\eta_i(\boldsymbol{\beta}) + r_i), Y_i) K_h(X_i - x_0),$$

where $\eta_i(\boldsymbol{\beta}) = \eta(X_i, \hat{\boldsymbol{\alpha}}) + (\mathbf{X}_i^T \boldsymbol{\beta} - \eta(x_0, \hat{\boldsymbol{\alpha}})) \frac{\eta(X_i, \hat{\boldsymbol{\alpha}})^\gamma}{\eta(x_0, \hat{\boldsymbol{\alpha}})^\gamma}$. The maximizer of the local quasi-log-likelihood $Q_u^*(\boldsymbol{\beta})$ is denoted by $\hat{\boldsymbol{\beta}}^* = \hat{\boldsymbol{\beta}}^*(x_0, \hat{\boldsymbol{\alpha}})$. Define

$$Q_u^{*\prime}(\boldsymbol{\beta}) \equiv \frac{\partial}{\partial \boldsymbol{\beta}} Q_u^*(\boldsymbol{\beta})$$

$$= \sum_{i=1}^n \frac{Y_i - g^{-1}(\eta_i(\boldsymbol{\beta}) + r_i)}{V(g^{-1}(\eta_i(\boldsymbol{\beta}) + r_i))} \times (g^{-1})'(\eta_i(\boldsymbol{\beta}) + r_i) \mathbf{X}_i K_h(X_i - x_0)$$

and similarly $Q_u^{*\prime\prime}(\boldsymbol{\beta}) = \frac{\partial^2}{\partial \boldsymbol{\beta} \partial \boldsymbol{\beta}^T} Q_u^*(\boldsymbol{\beta})$ to denote the gradient vector and Hessian matrix of the local quasi-likelihood $Q_u^*$, respectively. Applying Taylor's expansion to $Q_u^*(\boldsymbol{\beta})$ around $\hat{\boldsymbol{\beta}}(x_0, \hat{\boldsymbol{\alpha}})$, we get

$$0 = Q_u^{*\prime}(\hat{\boldsymbol{\beta}}^*) \approx Q_u^{*\prime}(\hat{\boldsymbol{\beta}}) + Q_u^{*\prime\prime}(\hat{\boldsymbol{\beta}})(\hat{\boldsymbol{\beta}}^* - \hat{\boldsymbol{\beta}})$$

which implies the following approximation of the estimation bias:

(6.1) $$\hat{\boldsymbol{\beta}}(x_0, \hat{\boldsymbol{\alpha}}) - \hat{\boldsymbol{\beta}}^*(x_0, \hat{\boldsymbol{\alpha}}) \approx (Q_u^{*\prime\prime}(\hat{\boldsymbol{\beta}}))^{-1} Q_u^{*\prime}(\hat{\boldsymbol{\beta}}).$$

Next we try to access the variance of the estimate $\hat{\boldsymbol{\beta}}$. To obtain variance, note that

$$0 = Q_u'(\hat{\boldsymbol{\beta}}) \approx Q_u'(\boldsymbol{\beta}^0) + Q_u''(\boldsymbol{\beta}^0)(\hat{\boldsymbol{\beta}} - \boldsymbol{\beta}^0),$$

where $\boldsymbol{\beta}^0 = (\beta_0^0, \beta_1^0, \ldots, \beta_p^0)^T$ with $\beta_j^0 = (\frac{\eta_0 - \eta(\cdot, \boldsymbol{\alpha})}{\eta(\cdot, \boldsymbol{\alpha})^\gamma})^{(p+j)}(x_0) \eta(X_i, \hat{\boldsymbol{\alpha}})^\gamma / j! + \eta(x_0, \boldsymbol{\alpha}) 1_{\{j=0\}}$ for $j = 0, 1, \ldots, p$. This implies that

$$\hat{\boldsymbol{\beta}} - \boldsymbol{\beta}^0 \approx -Q_u''(\boldsymbol{\beta}^0)^{-1} Q_u'(\boldsymbol{\beta}^0),$$

and an approximation for the conditional variance is given by

$$\text{var}(\hat{\boldsymbol{\beta}}|\mathbb{X}) \approx Q_u''(\boldsymbol{\beta}^0)^{-1} \text{var}(Q_u'(\boldsymbol{\beta}^0)|\mathbb{X}) Q_u''(\boldsymbol{\beta}^0)^{-1}.$$

Here the Hessian matrix can be approximated by $Q_u''(\hat{\boldsymbol{\beta}})$, and the variance term can be approximated as follows:

$$\text{var}(Q_u'(\boldsymbol{\beta}_0)|\mathbb{X}) = \sum_{i=1}^n \text{var}\left(\frac{\partial}{\partial \boldsymbol{\beta}} Q(g^{-1})(\eta_i(\boldsymbol{\beta}), Y_i) \Big| \mathbf{X}_i\right)_{\boldsymbol{\beta}=\boldsymbol{\beta}^0} K_h^2(X_i - x_0)$$

$$= \sum_{i=1}^n \xi_i \mathbf{X}_i \mathbf{X}_i^T K_h^2(X_i - x_0) \left(\frac{\eta(X_i, \hat{\boldsymbol{\alpha}})^\gamma}{\eta(x_0, \hat{\boldsymbol{\alpha}})^\gamma}\right)^2,$$

where

$$\xi_i = \operatorname{var}\left[\frac{Y_i - g^{-1}(\eta_i(\boldsymbol{\beta}))}{V(g^{-1}(\eta_i(\boldsymbol{\beta})))}(g^{-1})'(\eta_i(\boldsymbol{\beta}))\Big|\mathbf{X}_i\right]_{\boldsymbol{\beta} = \boldsymbol{\beta}^0}.$$

Note that $X_i$ has significant weight only in a neighborhood around $x_0$, and for such $i$,

$$\xi_i \approx [(g^{-1})'(\eta_0(x_0))]^2 / V(g^{-1}(\eta_0(x_0))).$$

Consequently, we have

$$(6.2) \qquad \operatorname{var}(Q_u'(\boldsymbol{\beta}_0)|\mathbb{X}) \approx \frac{[(g^{-1})'(\eta_0(x_0))]^2}{V(g^{-1}(\eta_0(x_0)))} S_n,$$

where

$$S_n = \sum_{i=1}^n \mathbf{X}_i \mathbf{X}_i^T K_h^2(X_i - x_0) \left(\frac{\eta(X_i, \hat{\boldsymbol{\alpha}})^\gamma}{\eta(x_0, \hat{\boldsymbol{\alpha}})^\gamma}\right)^2.$$

Combining the above results, we get

$$\operatorname{var}(\hat{\boldsymbol{\beta}}|\mathbb{X}) \approx \frac{[(g^{-1})'(\eta_0(x_0))]^2}{V(g^{-1}(\eta_0(x_0)))} Q_u''(\boldsymbol{\beta}^0)^{-1} S_n Q_u''(\boldsymbol{\beta}^0)^{-1},$$

where the unknown $\eta_0(x_0)$ and $\boldsymbol{\beta}^0$ can be replaced by their estimates $\hat{\eta}_u(x_0)$ and $\hat{\boldsymbol{\beta}}$, respectively.

6.2. *Bandwidth selection via bias-variance tradeoff.* Based on the above arguments, we first select a pilot bandwidth $\hat{h}_{p+a+1,p+a}^*$, which can be chosen using the ERSC. Next we fit a local polynomial with degree $p + a + 1$ and bandwidth $\hat{h}_{p+a+1,p+a}^*$ to get an estimate $\hat{\boldsymbol{\beta}}^{(p+a)} = (\hat{\beta}_0, \hat{\beta}_1, \ldots, \hat{\beta}_{p+a})^T$ via maximizing quasi-log-likelihood function (3.4). Using $\hat{\boldsymbol{\beta}}^{(p+a)}$, we get the approximation error $r_i$ and hence the estimated bias $\hat{B}_{p,0}(x; h)$ and variance $\hat{V}_{p,0}(x; h)$ of $\hat{\beta}_0$ which are respectively the first elements of the estimated bias vector (6.1) and variance matrix (6.2). An estimator of the mean squared error(MSE) of $\hat{\beta}_0$ is given by

$$\widehat{\mathrm{MSE}}_{p,0}(x_0; h) = \hat{B}_{p,0}^2(x_0; h) + \hat{V}_{p,0}(x_0; h)$$

which leads to our final bandwidth selector

$$(6.3) \qquad \hat{h}_{p,0} = \underset{h}{\operatorname{argmin}} \int \widehat{\mathrm{MSE}}_{p,0}(x; h)\, dx.$$

**7. Monte Carlo study.** In this section, we use simulations to illustrate the improvement of our newly proposed estimators by comparing them with the original nonparametric method. For simulations in this section and real data analysis in next section, we use the canonical link and local linear fitting by setting $p = 1$. To access the bias term in the pre-asymptotic bandwidth selection as discussed in Section 6, we choose the order of the approximation to the Taylor expansion error to be $a = 2$. In simulation studies, we first generate ten independent data sets, on which the pre-asymptotic bandwidth selector based on a grid search is applied. We set our final selected bandwidth to be the median of the obtained ten bandwidths, and it is fixed and used in our simulation. This speeds up the computation considerably. Different methods with their corresponding selected bandwidths are applied to another $R = 1000$ independent data sets and results are reported. When necessary, the Epanechnikov kernel is used in all of our numerical examples. For $\gamma$ in the unified family, we use a grid $\Gamma = 0, 0.1, 0.2, \ldots, 1, 1.2, 1.4, \ldots, 5$.

EXAMPLE 7.1 (Poisson). Each observation pair $(X, Y)$ in this example is generated in two steps: (1) the predictor variable $X$ is marginally uniformly distributed over $[-2, 2]$; (2) given $X = x$, the response $Y$ is generated from Poisson distribution with mean $\exp(\eta_0(x))$ where $\eta_0(x) = 3 \sin(\frac{\pi}{4}x - \frac{\pi}{2}) + 6$. Each sample consists of 100 i.i.d. pairs of observations. We estimate $\eta_0(\cdot)$ over $J = 100$ uniform grid points $\{x_j\}_{j=1}^J$ on the interval $[-2, 2]$. We use three different parametric guides: $G_1^P = \alpha_1 + \alpha_2 x + \alpha_3 x^2$, $G_2^P = \alpha_1 + \alpha_2 x + \alpha_3 x^2 + \alpha_4 x^3$ and $G_3^P = \alpha_1 + \alpha_2 \sin(\frac{\pi}{4}x - \frac{\pi}{2})$.

For an estimate $\hat{\eta}_r(\cdot)$ with $r$ indexing the replication of the simulation study, we define the bias $B_j = R^{-1} \sum_{r=1}^R [\hat{\eta}_r(x_j) - \eta_0(x_j)]$, the variance $S_j = R^{-1} \sum_{r=1}^R [\hat{\eta}_r(x_j) - R^{-1} \sum_{r'=1}^R \hat{\eta}_{r'}(x_j)]^2$ and the mean square error $\text{MSE}_j = B_j^2 + S_j$ at each $j$th grid point $x_j$. Let $B^2 = J^{-1} \sum_{j=1}^J B_j^2$, $V = J^{-1} \sum_{j=1}^J S_j$, $\text{MSE} = J^{-1} \sum_{j=1}^J \text{MSE}_j$ be the averages of the squared bias, variance and mean squared error (MSE) of the estimate $\hat{\eta}(\cdot)$, respectively. In Table 1, we report for different guides, the squared bias, variance and MSE of the original method, additive correction, multiplicative correction and the unified family of corrections with the best $\gamma$. The best $\gamma$ corresponds to the one that minimizes MSE over the grid $\Gamma$ and is given by 3.2, 3.0 and 3.0 for $G_1^P$, $G_2^P$ and $G_3^P$, respectively. The last block corresponds to $\gamma$ tuned by the selection method proposed in Section 5. The tuned $\gamma$ is given by 1.8, 1.8 and 0 for $G_1^P$, $G_2^P$ and $G_3^P$, respectively. The top panel "best $h$" means that each method uses its corresponding best smoothing bandwidth while the lower half "same $h$" corresponds to the case of using the same smoothing bandwidth selected by the original method.

The lower panel indicates that parametric guide reduces bias but has little effect on variance when the same $h$ is used for different methods. This is consistent with our asymptotic results. However, when the individual best $h$ is used for different methods, parametric guide also reduces variance as shown in the top half of Table 1. The underlying reason is that an appropriate parametric guide helps to make the nonparametric correction term flatter and smoother and a larger $h$ is allowed. This in turn reduces the variance. Table 1 indicates that the tuned $\gamma$ does not perform as well as the best $\gamma$, however, it improves over the additive and multiplicative corrections for either the quadratic or cubic guide. Note further that the improvement from the additive or multiplicative correction to the best $\gamma$ is much smaller than the improvement from the original method to the additive or multiplicative correction. Based on this observation, we recommend using the simple method outlined at the end of Section 5 to select $\gamma$ for real applications.

We plot $B^2$, $V$ and MSE for our parametrically-guided nonparametric estimation with different $\gamma$s in Figure 3 with the far-left isolated one corresponding to the original nonparametric estimation. Panels (A) and (C) correspond to cubic guide while panels (B) and (D) use the true sinusoid guide. The smoothing bandwidth is fixed at a same value for all estimation methods in panels (C) and (D) while each individual estimation uses its corresponding best smoothing parameter for panels (A) and (B). The figures give a picture of the results summarized in Table 1.

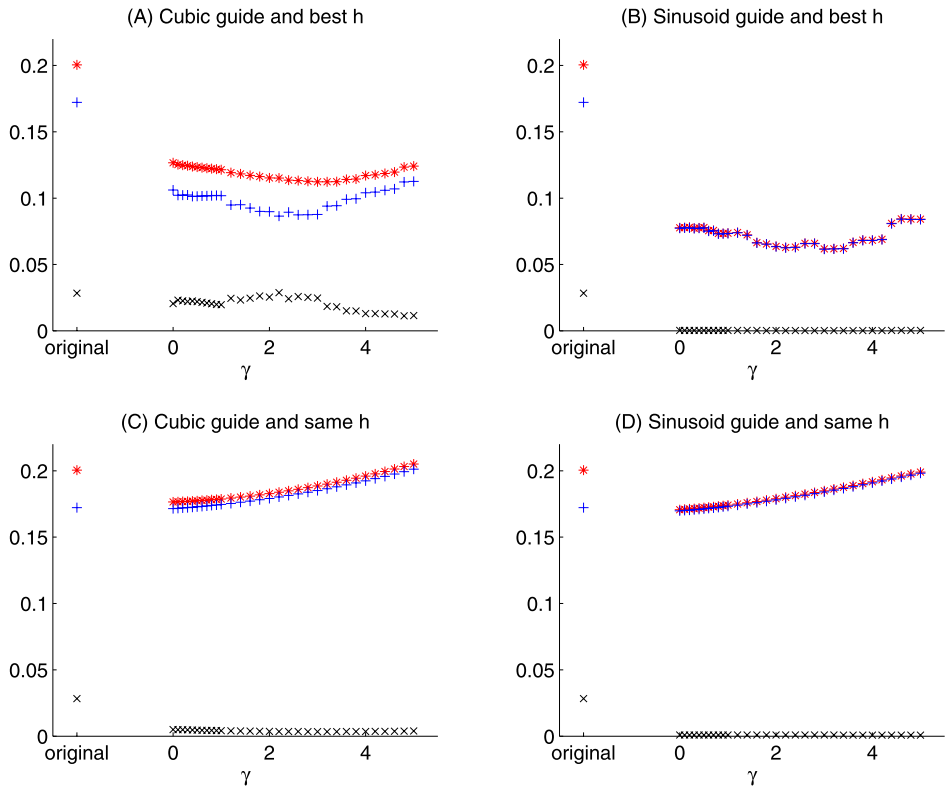

FIG. 3. Plots of $B^2$, $V$, and MSE denoted by black $\times$, blue $+$ and red $*$, respectively.

For a random sample of size 100, the best fit within the quadratic family is shown by the dotted line in panel (A) of Figure 1. The true unknown $\eta_0(\cdot)$, nonparametric estimate $\hat{\eta}(\cdot)$, two parametrically-guided estimates $\hat{\eta}_a(\cdot)$ and $\hat{\eta}_m(\cdot)$ are given by the solid, dotted, dashed and dot-dashed lines, respectively, in Figure 2. From this, we can see that parametrically-guided estimates improve the nonparametric counterpart around $x = 0$ where the curvature of $\eta_0(\cdot)$ is large and makes nonparametric estimation difficult.

EXAMPLE 7.2 (Bernoulli).   In this example, we consider Bernoulli distribution. The predictor variable $X$ is generated from Uniform$[-1, 1]$. Conditioning on $X = x$, the response $Y$ is generated from Bernoulli distribution with success probability $\exp(\eta_0(x))/(1 + \exp(\eta_0(x)))$ where $\eta_0(x) = 2\sin(\pi x)$. In this case, we consider samples of size 500 for two reasons: (1) the estimation of Bernoulli success probability is harder than the case of Poisson; (2) the use of a full sinusoid true $\eta_0(x)$ makes it even harder. Function $\eta_0(\cdot)$ is estimated over a uniform grid with $J = 100$ points over $[-1, 1]$. The average of squared bias, variance and MSE are reported in Table 2 for three guides $G_1^B = \alpha_1 + \alpha_2 x$, $G_2^B = \alpha_1 + \alpha_2 x + \alpha_3 x^2 + \alpha_4 x^3$ and $G_3^B = \alpha_1 + \alpha_2 \sin(\pi x)$. The best $\gamma$ is 1, 0.7 and 0.6 for $G_1^B$, $G_2^B$ and $G_3^B$, respectively. The tuned $\gamma$ is given by 0.8, 0.6 and 0.7 for $G_1^B$, $G_2^B$ and $G_3^B$, respectively.

Note that no improvement is observed for the additive correction with a linear guide $\alpha_1 + \alpha_2 x$ in Example 7.2. We can resort to our theoretical results to understand this exception. As we use local linear fitting, theoretically asymptotic bias depends on the second-order derivative of $\eta_0(\cdot) - \eta(\cdot, \boldsymbol{\alpha}_0) + \eta(x_0, \boldsymbol{\alpha}_0)$ and $\eta_0(\cdot)\eta(x_0, \boldsymbol{\alpha}_0)/\eta(\cdot, \boldsymbol{\alpha}_0)$ for additive and multiplicative corrections, respectively. A linear guide cannot reduce the second-order derivative of $\eta_0(\cdot) - \eta(\cdot, \boldsymbol{\alpha}_0) + \eta(x_0, \boldsymbol{\alpha}_0)$ and consequently does not reduce bias. However, a linear guide slightly reduces the second-order derivative of $\eta_0(\cdot)\eta(x_0, \boldsymbol{\alpha}_0)/\eta(\cdot, \boldsymbol{\alpha}_0)$ and improves the corresponding performance. This is consistent with our numerical results in Table 2. Note further that the multiplicative correction performs the best among the unified family of corrections when the linear guide is used.

**8. Real data analysis.**   In this section, we apply our newly proposed parametrically guided nonparametric estimation schemes to the Financial Aid Award Data, provided by National Longitudinal Survey of the High School Class of 1972. The data set is available online, and interested readers may find more information about this data set at http://www.oswego.edu/ kane/

TABLE 1
*Result on average squared bias, variance and MSE of Example 7.1*

| | Original method | | | | Additive $\gamma = 0$ | | | Multiplicative $\gamma = 1$ | | | Best $\gamma$ | | | Tuned $\gamma$ | | |
|---|---|---|---|---|---|---|---|---|---|---|---|---|---|---|---|---|
| | $B^2$ | $V$ | MSE | | $B^2$ | $V$ | MSE | $B^2$ | $V$ | MSE | $B^2$ | $V$ | MSE | $B^2$ | $V$ | MSE |
| Best $h$ | | | | $G_1^P$ | 2.05 | 10.24 | 12.28 | 1.97 | 9.84 | 11.81 | 2.12 | 8.87 | 10.99 | 3.18 | 8.14 | 11.32 |
| | 2.83 | 17.22 | 20.05 | $G_2^P$ | 2.05 | 10.61 | 12.66 | 1.96 | 10.18 | 12.14 | 2.46 | 8.77 | 11.23 | 2.87 | 8.64 | 11.51 |
| | | | | $G_3^P$ | 0.03 | 7.73 | 7.76 | 0.03 | 7.33 | 7.36 | 0.02 | 6.15 | 6.17 | 0.03 | 7.73 | 7.76 |
| Same $h$ | | | | $G_1^P$ | 0.50 | 17.00 | 17.50 | 0.42 | 17.33 | 17.74 | 0.50 | 17.00 | 17.50 | – | – | – |
| | 2.83 | 17.22 | 20.05 | $G_2^P$ | 0.50 | 17.15 | 17.65 | 0.41 | 17.46 | 17.87 | 0.50 | 17.15 | 17.65 | – | – | – |
| | | | | $G_3^P$ | 0.09 | 16.98 | 17.07 | 0.09 | 17.32 | 17.41 | 0.09 | 16.98 | 17.07 | – | – | – |

Note: All entries for squared bias, variance and MSE are multiplied by 100.

TABLE 2
*Result on the average of squared bias, variance, and MSE of Example 7.2*

| | Original method | | | | Additive $\gamma = 0$ | | | Multiplicative $\gamma = 1$ | | | Best $\gamma$ | | | Tuned $\gamma$ | | |
|---|---|---|---|---|---|---|---|---|---|---|---|---|---|---|---|---|
| | $B^2$ | $V$ | MSE | | $B^2$ | $V$ | MSE | $B^2$ | $V$ | MSE | $B^2$ | $V$ | MSE | $B^2$ | $V$ | MSE |
| | | | | $G_1^B$ | 92.4 | 721.6 | 814.0 | 88.1 | 719.7 | 807.7 | 88.1 | 719.7 | 807.7 | 88.1 | 719.7 | 807.7 |
| | 92.4 | 721.6 | 814.0 | $G_2^B$ | 160.5 | 510.7 | 671.1 | 49.2 | 636.1 | 685.4 | 91.4 | 551.6 | 643.0 | 94.9 | 548.6 | 643.5 |
| | | | | $G_3^B$ | 1.9 | 466.9 | 468.8 | 4.4 | 480.1 | 484.5 | 1.6 | 358.7 | 360.3 | 1.7 | 364.4 | 366.1 |

Note: All entries for squared bias, variance, and MSE are multiplied by 100.

econometrics/finaid.htm. There are twenty variables. We are interested in
using SAT score ($X$) to predict whether a student received financial aid
grants. There are 3076 students in total with SAT scores between 600 and
1300. Out of these 3076 students, 916 students received some financial aid
grants. The binary response $Y$ is coded in this way: $Y = 1$ means that a
student received financial aid grants and $Y = 0$ otherwise.

The pre-asymptotic bandwidth selector gives bandwidths 258.4615 for the
original nonparametric GLM. The corresponding nonparametric estimate of
the log odds ratio $\log \frac{P(Y=1|X=x)}{P(Y=0|X=x)}$ is given by the dot-dashed line in Figure
4. Based on this, we choose a cubic guide for two reasons. First, from the
simulation we know that the linear guide does not help at all in some cases.
Second, the nonparametric estimate does not indicate a quadratic shape.
Thus we apply the parametrically-guided logistic regression with a cubic
guide. The pre-asymptotic bandwidth selector gives bandwidths 296.1538
and 296.1538 for the parametrically-guided additive and parametrically-
guided multiplicative methods, respectively. Result is summarized in Figure
4. Cubic parametric estimate of the log odds ratio is given by the solid line;
our parametrically-guided estimates are given by the dashed and dotted lines
for additive and multiplicative methods, respectively.

We observe that our parametrically-guided additive and multiplicative es-
timates follow the cubic fit very closely. This suggests that there is no model
specification error by using a cubic model. However, the nonparametric es-
timate differs from the cubic fit for lower SAT scores.

**9. Discussion.** In this work, we extend the methodology of parametrically-
guided nonparametric estimation to GLMs and QLM. Asymptotic properties and numerical evidence demonstrate its improvement over the original non-parametric estimation scheme. There are possible extensions. For example, the whole estimation scheme can be easily extended to multivariate varying-coefficient and additive models. This enables us to incorporate prior knowledge into the analysis of multivariate nonparametric models, ameliorating the issues of curse of dimensionality.

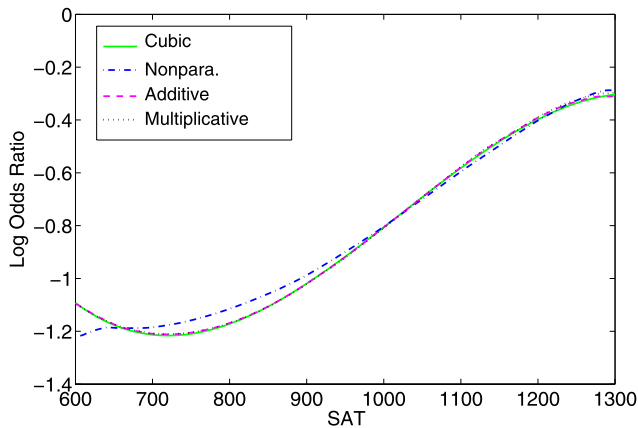

FIG. 4. *Plots of the nonparametric estimate, cubic estimate and our parametrical-
ly-guided estimates of the log odds ratio function for the Financial Aid Award Data.*

## APPENDIX: CONDITIONS AND PROOFS

Let $q_i(x, y) = (\partial^i / \partial x^i) Q(g^{-1}(x), y)$. Note that $q_i$ is linear in $y$ for fixed $x$ and that $q_1(\eta_0(x_0), \mu(x_0)) = 0$ and $q_2(\eta_0(x_0), \mu(x_0)) = -\rho(x_0)$.

The following technical conditions are imposed:

(A1) The function $q_2(x, y) < 0$ for $x \in \mathbb{R}$ and $y$ in the range of the response variable.

(A2) The functions $f'_X, \eta_0^{(p+2)}, \frac{\partial^{p+2}}{\partial x^{p+2}} \eta(x, \boldsymbol{\alpha}), \text{var}(Y|X = \cdot), V''$ and $g'''$ are continuous.

(A3) For each $x \in \text{supp}(f_X), \rho(x), \text{var}(Y|X = x)$ and $g'(\mu(x))$ are nonzero.

(A4) The kernel $K$ is a symmetric probability density with support $[-1, 1]$.

(A5) For each point $x_\delta$ on the boundary of $\text{supp}(f_X)$, there exists an interval $\mathcal{C}$ containing $x_\delta$ having nonnull interior such that $\inf_{x \in \mathcal{C}} f_X(x) > 0$.

White ([1982])-type conditions:

(B1) $E \log(f_0(x, y))$ exists and there exists a $m_1(x, y)$ such that $|\log(f(x, y; \boldsymbol{\alpha}))| \leq m_1(x, y)$ for any $\boldsymbol{\alpha} \in \mathbb{A}$ and $Em_1(x, y) < \infty$.

(B2) $E(\log(f_0(x, y)/f(x, y; \boldsymbol{\alpha})))$ has a unique minimizer $\boldsymbol{\alpha}_0$.

(B3) $\frac{\partial}{\partial \alpha_j} \log f(x, y; \boldsymbol{\alpha})$ is continuously differentiable in $\boldsymbol{\alpha}$ for $j = 1, 2, \ldots, q$.

(B4) There exist $m_2(x, y)$ and $m_3(x, y)$ such that $|\frac{\partial}{\partial \alpha_i} \log f(x, y; \boldsymbol{\alpha}) \frac{\partial}{\partial \alpha_j} \log f(x, y; \boldsymbol{\alpha})| \leq m_2(x, y)$ and $|\frac{\partial^2}{\partial \alpha_i \partial \alpha_j} \log f(x, y; \boldsymbol{\alpha})| \leq m_3(x, y)$ for any $\boldsymbol{\alpha} \in \mathbb{A}$, $1 \leq i, j \leq q$. Furthermore, both $Em_2(X, Y)$ and $Em_3(X, Y)$ exist.

(B5) Assume that $\boldsymbol{\alpha}_0$ is an interior point of $\mathbb{A}$; the matrix $(E \frac{\partial}{\partial \alpha_i} \log f(x, y; \boldsymbol{\alpha}) \frac{\partial}{\partial \alpha_j} \log f(x, y; \boldsymbol{\alpha}))_{1 \leq i, j \leq q}$ is nonsingular at $\boldsymbol{\alpha}_0$; $\boldsymbol{\alpha}_0$ is a regular point of matrix $(E \frac{\partial^2}{\partial \alpha_i \partial \alpha_j} \log f(x, y; \boldsymbol{\alpha}))_{1 \leq i, j \leq q}$.

For the case of unified correction with a fixed guide $\eta(x, \boldsymbol{\alpha})$, denote $\hat{\boldsymbol{\beta}} = \hat{\boldsymbol{\beta}}(x_0, \boldsymbol{\alpha}) = \text{argmax}_{\boldsymbol{\beta}} Q_u(\boldsymbol{\beta}; h, x_0, \boldsymbol{\alpha})$. Because $\hat{\boldsymbol{\beta}}$ is calculated using $X_i$ near $x_0$, we expect that

$$\eta(X_i, \hat{\boldsymbol{\alpha}}) + (\beta_0 + \cdots + \beta_p(X_i - x_0)^p - \eta(x_0, \hat{\boldsymbol{\alpha}})) \eta(X_i, \hat{\boldsymbol{\alpha}})^\gamma / \eta(x_0, \hat{\boldsymbol{\alpha}})^\gamma$$

$$\approx \eta_0(x_0) + \eta'_0(x_0)(X_i - x_0) + \cdots + \eta_0^{(p)}(x_0)(X_i - x_0)^p / p!.$$

Consequently, we expect that $\hat{\beta}_0 \to \eta_0(x_0)$ and

$$\hat{\beta}_j \to \left(\frac{\eta_0 - \eta(\cdot, \boldsymbol{\alpha})}{\eta(\cdot, \boldsymbol{\alpha})^\gamma}\right)^{(j)}(x_0)\eta(x_0, \boldsymbol{\alpha})^\gamma / j! \qquad \text{for } 1 \le j \le p.$$

We define $\phi_{\boldsymbol{\alpha},\gamma}(x) = (\eta_0(x) - \eta(x, \boldsymbol{\alpha}))/\eta(x, \boldsymbol{\alpha})^\gamma$ to simplify our notations. We thus study the asymptotic properties of

$$\hat{\boldsymbol{\beta}}^* = (nh)^{1/2}(\hat{\beta}_0 - \eta_0(x_0), h^1\{\hat{\beta}_1 - \eta(x_0, \boldsymbol{\alpha})^\gamma \phi_{\boldsymbol{\alpha},\gamma}^{(1)}(x_0)\}, \dots,$$
$$h^p\{p!\hat{\beta}_p - \eta(x_0, \boldsymbol{\alpha})^\gamma \phi_{\boldsymbol{\alpha},\gamma}^{(p)}(x_0)\})^T$$

so that each component has the same rate of convergence. Let $\mathbf{Q}_p(\mathcal{A})$ and $\mathbf{T}_p(\mathcal{A})$ be the $(p+1) \times (p+1)$ matrices having $(i,j)$th entry equal to $\nu_{i+j-1}(\mathcal{A})$ and $\int_{\mathcal{A}} z^{i+j-2} K^2(z) \, dz$. Also, define $\mathbf{D} = \text{diag}(1, 1/1!, \dots, 1/p!)$, $\boldsymbol{\Sigma}_x(\mathcal{A}) = \rho(x)f_X(x)\mathbf{D}\mathbf{N}_p(\mathcal{A})\mathbf{D}$,

$$\boldsymbol{\Gamma}_x(\mathcal{A}) = \frac{f_X(x)\text{var}(Y|X=x)}{\{V(\mu(x))g'(\mu(x))\}^2}\mathbf{D}\mathbf{T}_p(\mathcal{A})\mathbf{D},$$

$$\boldsymbol{\Lambda}_x(\mathcal{A}) = \frac{(\rho\eta^{2\gamma}(\cdot, \boldsymbol{\alpha})f_X)'(x)}{\eta^{2\gamma}(x, \boldsymbol{\alpha})}\mathbf{D}\mathbf{Q}_p(\mathcal{A})\mathbf{D},$$

$$a_{1,j}(\mathcal{A}) = \frac{1}{(p+1)!}\phi_{\boldsymbol{\alpha},\gamma}^{(p+1)}(x_0)\eta(x_0, \boldsymbol{\alpha})^\gamma \int_{\mathcal{A}} z^{p+1} K_{j-1,p}(z; \mathcal{A}) \, dz$$

and

$$a_{2,j}(\mathcal{A}) = \eta(x_0, \boldsymbol{\alpha})^\gamma \frac{1}{(p+2)!}\phi_{\boldsymbol{\alpha},\gamma}^{(p+2)}(x_0) \int_{\mathcal{A}} z^{p+2} K_{j-1,p}(z; \mathcal{A}) \, dz$$

$$+ \frac{\phi_{\boldsymbol{\alpha},\gamma}^{(p+1)}(x_0)(\rho\eta^{2\gamma}(\cdot, \boldsymbol{\alpha})f_X)'(x_0)}{(p+1)!(\rho\eta^\gamma(\cdot, \boldsymbol{\alpha})f_X)(x_0)}$$

$$\times \left\{\int_{\mathcal{A}} z^{p+2} K_{j-1,p}(z; \mathcal{A}) \, dz - (j-1) \int_{\mathcal{A}} z^{p+1} K_{j-2,p}(z; \mathcal{A}) \, dz\right.$$

$$\left. - \frac{1}{p!}\int_{\mathcal{A}} z^{p+1} K_{j-1,p}(z; \mathcal{A}) \, dz \int_{\mathcal{A}} z^{p+1} K_{p,p}(z; \mathcal{A}) \, dz\right\}.$$

Let $\mathbf{b}_{x_0}(\mathcal{A})$ be the $(p+1) \times 1$ vector having $j$th entry equal to $\sqrt{nh^{2p+3}}a_{1,j}(\mathcal{A}) + \sqrt{nh^{2p+5}}a_{2,j}(\mathcal{A})$.

MAIN THEOREM 1. *Suppose that conditions* (A1)–(A5) *hold and that* $h = h_n \to 0$, $nh^{2p+1} \to \infty$, $nh^{2p+3} < \infty$ *as* $n \to \infty$. *If* $x_0$ *is an interior point of* $\text{supp}(f_X)$, *and* $p > 0$, *then*

$$\{\boldsymbol{\Sigma}_{x_0}([-1,1])^{-1}\boldsymbol{\Gamma}_{x_0}([-1,1])\boldsymbol{\Sigma}_{x_0}([-1,1])^{-1}\}^{-1/2}$$
$$\times \{\hat{\boldsymbol{\beta}}^* - \mathbf{b}_{x_0}([-1,1]) + o(\sqrt{nh^{2p+5}})\} \xrightarrow{\mathcal{D}} N(\mathbf{0}, \mathbf{I}_{p+1}).$$

*If* $x_0 = x_n$ *is of the form* $x_0 = x_\delta + hc$ *where* $c \in [-1, 1]$ *is fixed and* $x_\delta$ *is a fixed point on the boundary of* $\operatorname{supp}(f_X)$, *then*

$$\{\boldsymbol{\Sigma}_{x_0}(\mathcal{D}_{x_0,h})^{-1}\boldsymbol{\Gamma}_{x_0}(\mathcal{D}_{x_0,h})\boldsymbol{\Sigma}_{x_0}(\mathcal{D}_{x_0,h})^{-1}\}^{-1/2}$$
$$\times\{\hat{\boldsymbol{\beta}}^* - \mathbf{b}_{x_0}(\mathcal{D}_{x_0,h}) + o(\sqrt{nh^{2p+5}})\} \stackrel{\mathcal{D}}{\to} N(\mathbf{0}, \mathbf{I}_{p+1}).$$

The proof of the main theorem follows directly from Lemmas 1 and 2, which are stated and proved as follows. Denote $\mathbf{Z}_i = (1, (X_i - x_0)/h, \ldots, (X_i - x_0)^p/(h^p p!))^T$.

LEMMA 1. *Let* $\bar{\eta}(x_0, x) = \eta(x_0, \boldsymbol{\alpha})^\gamma \sum_{j=0}^p \phi_{\boldsymbol{\alpha},\gamma}^{(j)}(x_0)(x - x_0)^j/j!$ *and* $\mathbf{W}_n = (nh)^{-1/2} \sum_{i=1}^n \mathbf{Y}_i^*$ *where*

$$\mathbf{Y}_i^* = \left(\frac{\eta(X_i, \boldsymbol{\alpha})}{\eta(x_0, \boldsymbol{\alpha})}\right)^\gamma q_1\left(\eta(X_i, \boldsymbol{\alpha}) + \left(\frac{\eta(X_i, \boldsymbol{\alpha})}{\eta(x_0, \boldsymbol{\alpha})}\right)^\gamma \bar{\eta}(x_0, X_i), Y_i\right) K\left(\frac{X_i - x_0}{h}\right) \mathbf{Z}_i.$$

*Then under conditions* (A1)–(A5), $nh^3 \to \infty$ *and* $h \to 0$, *we have*

$$\hat{\boldsymbol{\beta}}^* = \boldsymbol{\Sigma}_{x_0}^{-1} \mathbf{W}_n - h\boldsymbol{\Sigma}_{x_0}^{-1} \boldsymbol{\Lambda}_{x_0} \boldsymbol{\Sigma}_{x_0}^{-1} \mathbf{W}_n + o_P(h).$$

PROOF. Recall that $\hat{\boldsymbol{\beta}}$ maximizes $Q_u(\boldsymbol{\beta}; x_0, p, \boldsymbol{\alpha})$. Let

$$\boldsymbol{\beta}^* \equiv (nh)^{1/2}(\beta_0 - \eta(x_0), h^1\{\beta_1 - \eta(x_0, \boldsymbol{\alpha})^\gamma \phi_{\boldsymbol{\alpha},\gamma}^{(1)}(x_0)\},$$
$$\ldots, h^p\{p!\beta_p - \eta(x_0, \boldsymbol{\alpha})^\gamma \phi_{\boldsymbol{\alpha},\gamma}^{(p)}(x_0)\})^T.$$

Then

$$\eta(X_i, \boldsymbol{\alpha}) + \frac{\eta(X_i, \boldsymbol{\alpha})^\gamma}{\eta(x_0, \boldsymbol{\alpha})^\gamma}\{\beta_0 + \beta_1(X_i - x_0) + \cdots + \beta_p(X_i - x_0)^p - \eta(x_0, \hat{\boldsymbol{\alpha}})\}$$

$$= \eta(X_i, \boldsymbol{\alpha}) + \frac{\eta(X_i, \boldsymbol{\alpha})^\gamma}{\eta(x_0, \boldsymbol{\alpha})^\gamma}\{\bar{\eta}(x_0, X_i) + a_n \boldsymbol{\beta}^{*T} \mathbf{Z}_i\},$$

where $a_n = (nh)^{-1/2}$. If $\hat{\boldsymbol{\beta}}$ maximizes $Q_u(\boldsymbol{\beta}; x_0, p, \boldsymbol{\alpha})$, then $\hat{\boldsymbol{\beta}}^*$ maximizes

$$\sum_{i=1}^n Q\left(g^{-1}\left(\eta(X_i, \boldsymbol{\alpha}) + \frac{\eta(X_i, \boldsymbol{\alpha})^\gamma}{\eta(x_0, \boldsymbol{\alpha})^\gamma}\{\bar{\eta}(x_0, X_i) + a_n \boldsymbol{\beta}^{*T} \mathbf{Z}_i\}\right), Y_i\right) K\left(\frac{X_i - x_0}{h}\right)$$

as a function of $\boldsymbol{\beta}^*$. To study the asymptotic properties of $\hat{\boldsymbol{\beta}}^*$, we apply the quadratic approximation lemma [Fan and Gijbels (1995)] to the maximization of the normalized function

$$l_n(\boldsymbol{\beta}^*) = \sum_{i=1}^n \left\{Q\left(g^{-1}\left(\eta(X_i, \boldsymbol{\alpha}) + \left(\frac{\eta(X_i, \boldsymbol{\alpha})}{\eta(x_0, \boldsymbol{\alpha})}\right)^\gamma (\bar{\eta}(x_0, X_i) + a_n \boldsymbol{\beta}^{*T} \mathbf{Z}_i)\right), Y_i\right)\right.$$

$$-Q\left(g^{-1}\left(\eta(X_i,\boldsymbol{\alpha})+\left(\frac{\eta(X_i,\boldsymbol{\alpha})}{\eta(x_0,\boldsymbol{\alpha})}\right)^\gamma\bar\eta(x_0,X_i)\right),Y_i\right)\right\}$$

$$\times K\left(\frac{X_i-x_0}{h}\right).$$

Then $\hat{\boldsymbol{\beta}}^*$ maximize $l_n$. We remark that condition (A1) implies that $l_n$ is concave in $\boldsymbol{\beta}^*$. Using a Taylor series expansion of $Q(g^{-1}(\cdot),Y_i)$,

$$
\begin{aligned}
l_n(\boldsymbol{\beta}^*) &= a_n\sum_{i=1}^n\left(\frac{\eta(X_i,\boldsymbol{\alpha})}{\eta(x_0,\boldsymbol{\alpha})}\right)^\gamma q_1\left(\eta(X_i,\boldsymbol{\alpha})+\left(\frac{\eta(X_i,\boldsymbol{\alpha})}{\eta(x_0,\boldsymbol{\alpha})}\right)^\gamma\bar\eta(x_0,X_i),Y_i\right)\\
&\qquad\times\boldsymbol{\beta}^{*T}\mathbf{Z}_iK\{(X_i-x_0)/h\}\\
&\quad+\frac{a_n^2}{2}\sum_{i=1}^n\left(\frac{\eta(X_i,\boldsymbol{\alpha})}{\eta(x_0,\boldsymbol{\alpha})}\right)^{2\gamma}q_2\left(\eta(X_i,\boldsymbol{\alpha})+\left(\frac{\eta(X_i,\boldsymbol{\alpha})}{\eta(x_0,\boldsymbol{\alpha})}\right)^\gamma\bar\eta(x_0,X_i),Y_i\right)\\
&\qquad\times(\boldsymbol{\beta}^{*T}\mathbf{Z}_i)^2K\{(X_i-x_0)/h\}\\
&\quad+\frac{a_n^3}{6}\sum_{i=1}^n\left(\frac{\eta(X_i,\boldsymbol{\alpha})}{\eta(x_0,\boldsymbol{\alpha})}\right)^{3\gamma}q_3\left(\eta(X_i,\boldsymbol{\alpha})+\left(\frac{\eta(X_i,\boldsymbol{\alpha})}{\eta(x_0,\boldsymbol{\alpha})}\right)^\gamma\eta_i,Y_i\right)(\boldsymbol{\beta}^{*T}\mathbf{Z}_i)^3\\
&\qquad\times K\{(X_i-x_0)/h\},
\end{aligned}
$$
(A.1)

where $\eta_i$ is between $\bar\eta(x_0,X_i)$ and $\bar\eta(x_0,X_i)+a_n\boldsymbol{\beta}^{*T}\mathbf{Z}_i$. Let

$$
\begin{aligned}
\mathbf{A}_n &= a_n^2\sum_{i=1}^n\left(\frac{\eta(X_i,\boldsymbol{\alpha})}{\eta(x_0,\boldsymbol{\alpha})}\right)^{2\gamma}q_2\left(\eta(X_i,\boldsymbol{\alpha})+\left(\frac{\eta(X_i,\boldsymbol{\alpha})}{\eta(x_0,\boldsymbol{\alpha})}\right)^\gamma\bar\eta(x_0,X_i),Y_i\right)\\
&\qquad\times K\{(X_i-x_0)/h\}\mathbf{Z}_i\mathbf{Z}_i^T.
\end{aligned}
$$

Then the second term in (A.1) equals $\frac{1}{2}\boldsymbol{\beta}^{*T}\mathbf{A}_n\boldsymbol{\beta}^*$. Now $(\mathbf{A}_n)_{ij}=(E\mathbf{A}_n)_{ij}+O_P(\sqrt{\mathrm{var}((\mathbf{A}_n)_{ij})})$ and

$$
\begin{aligned}
E\mathbf{A}_n &= h^{-1}E\left\{\left(\frac{\eta(X_1,\boldsymbol{\alpha})}{\eta(x_0,\boldsymbol{\alpha})}\right)^{2\gamma}q_2\left(\eta(X_1,\boldsymbol{\alpha})+\left(\frac{\eta(X_1,\boldsymbol{\alpha})}{\eta(x_0,\boldsymbol{\alpha})}\right)^\gamma\bar\eta(x_0,X_1),\mu(X_1)\right)\right.\\
&\qquad\qquad\left.\times K\{(X_1-x_0)/h\}\mathbf{Z}_1\mathbf{Z}_1^T\right\}
\end{aligned}
$$

since $q_2$ is linear in $y$ for fixed $x$. Because $\mathrm{supp}(K)=[-1,1]$, we need only consider $|X_1-x_0|\le h$, and thus

$$
\begin{aligned}
&\eta(X_1,\boldsymbol{\alpha})+\left(\frac{\eta(X_1,\boldsymbol{\alpha})}{\eta(x_0,\boldsymbol{\alpha})}\right)^\gamma\bar\eta(x_0,X_1)-\eta_0(X_1)\\
&\quad=-\eta(X_1,\boldsymbol{\alpha})^\gamma\left\{\frac{1}{(p+1)!}\phi_{\boldsymbol{\alpha},\gamma}^{(p+1)}(x_0)(X_1-x_0)^{p+1}\right.
\end{aligned}
$$

$$+ \frac{1}{(p+2)!} \phi_{\boldsymbol{\alpha},\gamma}^{(p+2)}(x_0)(X_1 - x_0)^{p+2} \Big\} + o(h^{p+2}).$$

Then

$$(i-1)!(j-1)!(E\mathbf{A}_n)_{ij}$$

$$= h^{-1} E \Bigg\{ \left( \frac{\eta(X_1, \boldsymbol{\alpha})}{\eta(x_0, \boldsymbol{\alpha})} \right)^{2\gamma}$$

$$\times q_2 \Big( \eta(X_1, \boldsymbol{\alpha}) + \left( \frac{\eta(X_1, \boldsymbol{\alpha})}{\eta(x_0, \boldsymbol{\alpha})} \right)^{\gamma} \bar{\eta}(x_0, X_1), \mu(X_1) \Big)$$

$$\times K \left( \frac{X_1 - x_0}{h} \right) \left( \frac{X_1 - x_0}{h} \right)^{i+j-2} \Bigg\}$$

$$= \int \left( \frac{\eta(x_0 + hZ, \boldsymbol{\alpha})}{\eta(x_0, \boldsymbol{\alpha})} \right)^{2\gamma}$$

$$\times q_2 \Big( \eta(x_0 + hZ, \boldsymbol{\alpha})$$

$$+ \left( \frac{\eta(x_0 + hZ, \boldsymbol{\alpha})}{\eta(x_0, \boldsymbol{\alpha})} \right)^{\gamma} \bar{\eta}(x_0, x_0 + hZ), \mu(x_0 + hZ) \Big)$$

$$\times K(Z) Z^{i+j-2} f_X(x_0 + hZ) \, dZ$$

$$= \int \left( \frac{\eta(x_0 + hZ, \boldsymbol{\alpha})}{\eta(x_0, \boldsymbol{\alpha})} \right)^{2\gamma} q_2(\eta_0(x_0 + hZ) + o(h^p), \mu(x_0 + hZ))$$

$$\times K(Z) Z^{i+j-2} f_X(x_0 + hZ) \, dZ$$

$$= \int \left( \frac{\eta(x_0 + hZ, \boldsymbol{\alpha})}{\eta(x_0, \boldsymbol{\alpha})} \right)^{2\gamma} [q_2(\eta_0(x_0 + hZ), \mu(x_0 + hZ)) + o(h^p)]$$

$$\times K(Z) Z^{i+j-2} f_X(x_0 + hZ) \, dZ$$

$$= \int - \left( \frac{\eta(x_0 + hZ, \boldsymbol{\alpha})}{\eta(x_0, \boldsymbol{\alpha})} \right)^{2\gamma} \rho(x_0 + hZ) f_X(x_0 + hZ)$$

$$\times K(Z) Z^{i+j-2} \, dZ + o(h)$$

$$= -(\rho f_X)(x_0) \nu_{i+j-2} - h \frac{(\rho \eta^{2\gamma}(\cdot, \boldsymbol{\alpha}) f_X)'(x_0)}{\eta^{2\gamma}(x_0, \boldsymbol{\alpha})} \nu_{i+j-1} + o(h).$$

Similar arguments show that $\text{var}\{(\mathbf{A}_n)_{ij}\} = O\{(nh)^{-1}\}$ and that the last term in (A.1) is $O_P\{(nh)^{-1/2}\}$. Therefore, $l_n(\boldsymbol{\beta}^*) = \mathbf{W}_n^T \boldsymbol{\beta}^* - \frac{1}{2} \boldsymbol{\beta}^{*T} (\boldsymbol{\Sigma}_{x_0} + h\boldsymbol{\Lambda}_{x_0}) \boldsymbol{\beta}^* + o_P(h)$ because $nh^3 \to \infty$ and $h \to 0$. Similar arguments show that $l_n'(\boldsymbol{\beta}^*) = \mathbf{W}_n - (\boldsymbol{\Sigma}_{x_0} + h\boldsymbol{\Lambda}_{x_0}) \boldsymbol{\beta}^* + o_P(h)$ and $l_n''(\boldsymbol{\beta}^*) = -(\boldsymbol{\Sigma}_{x_0} + h\boldsymbol{\Lambda}_{x_0}) + o_P(h)$. The result follows directly from the quadratic approximation lemma of Fan and Gijbels (1995). $\square$

Lemma 2. *Suppose that the conditions of Theorem 1 hold. For $\mathbf{W}_n$ as defined in Lemma 1,*

$$\{\mathbf{\Sigma}_{x_0}^{-1} - h\mathbf{\Sigma}_{x_0}^{-1}\mathbf{\Lambda}_{x_0}\mathbf{\Sigma}_{x_0}^{-1}\}E(\mathbf{W}_n) = \mathbf{b}_{x_0} + o\{(nh^{2p+5})^{1/2}\},$$

$$\mathbf{\Gamma}_{x_0}^{-1/2}\operatorname{cov}(\mathbf{W}_n)\mathbf{\Gamma}_{x_0}^{-1/2} \to \mathbf{I}_{p+1}$$

*and*

$$\mathbf{\Gamma}_{x_0}^{-1/2}(\mathbf{W}_n - E\mathbf{W}_n) \xrightarrow{\mathcal{D}} N(\mathbf{0}, \mathbf{I}_{p+1}).$$

Proof. We compute the mean and covariance matrix of the random vector $\mathbf{W}_n$ by studying $\mathbf{Y}_1^*$, as defined in Lemma 1. Denote $(E\mathbf{Y}_1^*)_i$ to be the mean of the $i$th component of $\mathbf{Y}_1^*$. Then it is easy to show that $\frac{(i-1)!}{h}(E\mathbf{Y}_1^*)_i$ is equal to

$$\int \left(\frac{\eta(x_0 + hZ, \boldsymbol{\alpha})}{\eta(x_0, \boldsymbol{\alpha})}\right)^\gamma q_1\bigg(\eta(x_0 + hZ, \boldsymbol{\alpha})$$
$$+ \left(\frac{\eta(x_0 + hZ, \boldsymbol{\alpha})}{\eta(x_0, \boldsymbol{\alpha})}\right)^\gamma \bar{\eta}(x_0, x_0 + hZ), \mu(x_0 + hZ)\bigg)$$
$$\times Z^{i-1}K(Z)f_X(x_0 + hZ)\, dZ.$$

Now by the Taylor expansion,

$$q_1\bigg(\eta(x_0 + hZ, \boldsymbol{\alpha}) + \left(\frac{\eta(x_0 + hZ, \boldsymbol{\alpha})}{\eta(x_0, \boldsymbol{\alpha})}\right)^\gamma \bar{\eta}(x_0, x_0 + hZ), \mu(x_0 + hZ)\bigg)$$
$$= \eta(x_0 + hZ, \boldsymbol{\alpha})^\gamma \bigg\{\frac{\phi_{\boldsymbol{\alpha}, \gamma}^{(p+1)}(x_0)}{(p+1)!}(hZ)^{p+1}$$
$$+ \frac{\phi_{\boldsymbol{\alpha}, \gamma}^{(p+2)}(x_0)}{(p+2)!}(hZ)^{p+2} + o(h^{p+2})\bigg\}\rho(x_0 + hZ)$$
$$+ o(h^{p+2}).$$

Thus

$$(E\mathbf{Y}_1^*)_i = \frac{(\rho\eta^{2\gamma}(\cdot, \boldsymbol{\alpha})f_X)(x_0)}{\eta(x_0, \boldsymbol{\alpha})^\gamma(i-1)!}\left(h^{p+2}\frac{\phi_{\boldsymbol{\alpha}, \gamma}^{(p+1)}(x_0)}{(p+1)!}\nu_{p+i} + h^{p+3}\zeta_p(x_0)\nu_{p+i+1}\right)$$
$$\text{(A.2)} \qquad + o(h^{p+3}),$$

where

$$\zeta_p(x_0) = \frac{1}{(p+2)!}\phi_{\boldsymbol{\alpha}, \gamma}^{(p+2)}(x_0) + \frac{\phi_{\boldsymbol{\alpha}, \gamma}^{(p+1)}(x_0)(\rho\eta^{2\gamma}(\cdot, \boldsymbol{\alpha})f_X)'(x_0)}{(p+1)!(\rho\eta^{2\gamma}(\cdot, \boldsymbol{\alpha})f_X)(x_0)}.$$

Note that

$$
\begin{aligned}
\boldsymbol{\Sigma}_{x_0}^{-1} E \mathbf{W}_n &= \left(\frac{n}{h}\right)^{1/2} \frac{1}{(\rho f_X)(x_0)} D^{-1} N^{-1} D^{-1} E \mathbf{Y}_1^* \\
&= \left(\frac{n}{h}\right)^{1/2} D^{-1} N^{-1} \Big[ h^{p+2} \frac{1}{(p+1)!} \phi_{\boldsymbol{\alpha},\gamma}^{(p+1)}(x_0) \eta(x_0,\boldsymbol{\alpha})^\gamma \\
&\qquad\qquad \times (\nu_{p+1}, \nu_{p+2}, \ldots, \nu_{2p+1})^T \\
&\qquad\qquad + h^{p+3} \zeta_p(x_0) \eta(x_0,\boldsymbol{\alpha})^\gamma (\nu_{p+2}, \nu_{p+3}, \ldots, \nu_{2p+2})^T \Big].
\end{aligned}
$$

The $i$th component of $\boldsymbol{\Sigma}_{x_0}^{-1} E \mathbf{W}_n$ is

$$
\begin{aligned}
(\boldsymbol{\Sigma}_{x_0}^{-1} E \mathbf{W}_n)_i &= \left(\frac{n}{h}\right)^{1/2} (i-1)! (\{N^{-1}\}_{i,1}, \{N^{-1}\}_{i,2}, \ldots, \{N^{-1}\}_{i,p+1}) \\
&\qquad \times \Big[ h^{p+2} \frac{1}{(p+1)!} \phi_{\boldsymbol{\alpha},\gamma}^{(p+1)}(x_0) \eta(x_0,\boldsymbol{\alpha})^\gamma (\nu_{p+1}, \nu_{p+2}, \ldots, \nu_{2p+1})^T \\
&\qquad\qquad + h^{p+3} \zeta_p(x_0) \eta(x_0,\boldsymbol{\alpha})^\gamma (\nu_{p+2}, \nu_{p+3}, \ldots, \nu_{2p+2})^T \Big] \\
&= (nh^{2p+3})^{1/2} \frac{1}{(p+1)!} \phi_{\boldsymbol{\alpha},\gamma}^{(p+1)}(x_0) \eta(x_0,\boldsymbol{\alpha})^\gamma \int z^{p+1} K_{i-1,p}(z)\, dz \\
&\quad + (nh^{2p+5})^{1/2} \zeta_p(x_0) \eta(x_0,\boldsymbol{\alpha})^\gamma \int z^{p+2} K_{i-1,p}(z)\, dz \\
&\quad + o\{(nh^{2p+5})^{1/2}\}.
\end{aligned}
$$

Next, consider the second term in the expression

$$
\begin{aligned}
h(\boldsymbol{\Sigma}_{x_0}^{-1} \boldsymbol{\Lambda}_{x_0} \boldsymbol{\Sigma}_{x_0}^{-1} E \mathbf{W}_n)_i &= (nh^{2p+5})^{1/2} (i-1)! \frac{\phi_{\boldsymbol{\alpha},\gamma}^{(p+1)}(x_0)(\rho \eta^{2\gamma}(\cdot,\boldsymbol{\alpha}) f_X)'(x_0)}{(p+1)!(\rho \eta^\gamma(\cdot,\boldsymbol{\alpha}) f_X)(x_0)} \\
&\qquad \times \sum_{j=1}^{p+1} (\mathbf{N}_p^{-1} \mathbf{Q}_p \mathbf{N}_p^{-1})_{ij} \nu_{p+j} + O\{(nh^{2p+7})^{1/2}\}.
\end{aligned}
$$

Using the fact that $(\mathbf{Q}_p)_{kl} = (\mathbf{N}_p)_{k,l+1}$ for $l < p+1$, it can be shown that for $i = 2, \ldots, p+1$,

$$
(\mathbf{N}_p^{-1} \mathbf{Q}_p \mathbf{N}_p^{-1})_{ij} = (\mathbf{N}_p^{-1})_{i-1,j} + \left\{ \sum_{k=1}^{p+1} (\mathbf{N}_p^{-1})_{i,k} \nu_{p+k} \right\} (\mathbf{N}_p^{-1})_{p+1,j}
$$

and by similar reasoning,

$$
(\mathbf{N}_p^{-1} \mathbf{Q}_p \mathbf{N}_p^{-1})_{1j} = \left\{ \sum_{k=1}^{p+1} (\mathbf{N}_p^{-1})_{1,k} \nu_{p+k} \right\} (\mathbf{N}_p^{-1})_{p+1,j}.
$$

So by Lemma 3 of Fan, Heckman and Wand (1995),

$$(i-1)! \sum_{j=1}^{p+1} (\mathbf{N}_p^{-1} \mathbf{Q}_p \mathbf{N}_p^{-1})_{ij} \nu_{p+j}$$

$$= (i-1) \int z^{p+1} K_{i-2,p}(z) \, dz$$

$$+ \frac{1}{p!} \int z^{p+1} K_{p,p}(z) \, dz \int z^{p+1} K_{i-1,p}(z) \, dz.$$

The statement concerning the asymptotic mean follows immediately. By (A.2), the covariance between the $i$th and $j$th component of $\mathbf{Y}_1^*$ is $E((\mathbf{Y}_1^*)_i (\mathbf{Y}_1^*)_j) + O(h^{2p+4})$. By a Taylor series expansion, $E((\mathbf{Y}_1^*)_i (\mathbf{Y}_1^*)_j)$ is given by

$$E\left[ \left( \left( \frac{\eta(X_1, \boldsymbol{\alpha})}{\eta(x_0, \boldsymbol{\alpha})} \right)^\gamma q_1 \left( \eta(X_1, \boldsymbol{\alpha}) + \left( \frac{\eta(X_1, \boldsymbol{\alpha})}{\eta(x_0, \boldsymbol{\alpha})} \right)^\gamma \bar{\eta}(x_0, X_1), Y_1 \right) K\{(X_1 - x_0)/h\} \right)^2 \right.$$

$$\left. \times \frac{\{(X_1 - x_0)/h\}^{i+j-2}}{(i-1)!(j-1)!} \right]$$

$$= \int \left( \left( \frac{\eta(x_0 + hZ, \boldsymbol{\alpha})}{\eta(x_0, \boldsymbol{\alpha})} \right)^\gamma q_1 \left( \eta(x_0 + hZ, \boldsymbol{\alpha}) + \left( \frac{\eta(x_0 + hZ, \boldsymbol{\alpha})}{\eta(x_0, \boldsymbol{\alpha})} \right)^\gamma \right.\right.$$

$$\left.\left. \times \bar{\eta}(x_0, x_0 + hZ), Y_1 \right) K(Z) \right)^2$$

$$\times \frac{Z^{i+j-2} f_X(x_0 + hZ) h}{(i-1)!(j-1)!} \, dZ$$

$$= \int (q_1(\eta(x_0 + hZ), Y_1) K(Z))^2 \frac{Z^{i+j-2} f_X(x_0) h}{(i-1)!(j-1)!} \, dZ + o(h).$$

Noticing that

$$q_1(\eta(x_0 + hZ), Y_1)$$

$$= \frac{Y_1 - g^{-1}(\eta(x_0 + hZ))}{V(g^{-1}(\eta(x_0 + hZ)))} (g^{-1})'(\eta(x_0 + hz)),$$

we can derive

$$\{\text{cov}(\mathbf{Y}_1^*)\}_{ij} = \frac{h f_X(x_0) \text{var}(Y | X = x_0)}{[V(\mu(x_0)) g'(\mu(x_0))]^2} \int \frac{z^{i+j-2}}{(i-1)!(j-1)!} K^2(Z) \, dZ + o(h).$$

Therefore, $\boldsymbol{\Gamma}_{x_0}^{-1/2} \text{cov}(\mathbf{W}_n) \boldsymbol{\Gamma}_{x_0}^{-1/2} \to \mathbf{I}_{p+1}$. Now, we use the Cramér–Wold device to derive the asymptotic normality of $\mathbf{W}_n$. For any unit vector $\mathbf{u} \in \mathbb{R}^{p+1}$, if

$$(A.3) \qquad (na_n^2)^{-1/2} \mathbf{u}^T \text{cov}(\mathbf{Y}_1^*)^{-1/2} (\mathbf{W}_n - E\mathbf{W}_n) \to_D N(0, 1)$$

then $h^{1/2} \operatorname{cov}(\mathbf{Y}_1^*)^{-1/2}(\mathbf{W}_n - E\mathbf{W}_n) \to_D N(\mathbf{0}, \mathbf{I}_{p+1})$, and so $\mathbf{\Gamma}_{x_0}^{-1/2}(\mathbf{W}_n - E\mathbf{W}_n) \to_D N(\mathbf{0}, \mathbf{I}_{p+1})$. To prove (A.3), we only need to check Lyapounov's condition for that sequence which can be easily verified. $\square$

Noting that $\hat{\beta}_0 \to \eta_0(x_0)$ and $\hat{\beta}_j \to \eta(x_0, \boldsymbol{\alpha})^\gamma (\frac{\eta_0 - \eta(\cdot, \boldsymbol{\alpha})}{\eta(\cdot, \boldsymbol{\alpha})^\gamma})^{(j)}(x_0)/j!$ for $1 \le j \le p$, we can define the estimator of $\eta_0^{(j)}(x_0)$ iteratively as $\hat{\eta}_{u,0}(x_0; p, h, \boldsymbol{\alpha}) \equiv \hat{\eta}_u(x_0; p, h, \boldsymbol{\alpha}) = \hat{\beta}_0$ and

$$\hat{\eta}_{u,j}(x_0; p, h) = j!\hat{\beta}_j - \eta(x_0, \boldsymbol{\alpha})^\gamma \sum_{i=0}^{j-1} \hat{\eta}_{u,i}(x_0; p, h)(1/\eta^\gamma)^{(j-i)}(x_0, \boldsymbol{\alpha}) \binom{j}{i}$$
$$+ \eta(x_0, \boldsymbol{\alpha})^\gamma (\eta^{1-\gamma})^{(j)}(x_0, \boldsymbol{\alpha}) \qquad \text{for } 1 \le j \le p,$$

where $\binom{j}{i} = j!/(i!(j-i)!)$. Simple algebra leads to

$$\hat{\eta}_{u,j}(x_0; p, h) - \eta_0^{(j)}(x_0)$$
$$= j!\hat{\beta}_j - \eta(x_0, \boldsymbol{\alpha})^\gamma \left(\frac{\eta_0 - \eta(\cdot, \boldsymbol{\alpha})}{\eta(\cdot, \boldsymbol{\alpha})^\gamma}\right)^{(j)}(x_0)$$
$$- \eta(x_0, \boldsymbol{\alpha})^\gamma \sum_{i=0}^{j-1} (\hat{\eta}_{u,i}(x_0; p, h) - \eta_0^{(i)}(x_0))$$
$$\times (1/\eta^\gamma)^{(j-i)}(x_0, \boldsymbol{\alpha}) \binom{j}{i} \qquad \text{for } 1 \le j \le p.$$

Denote $w_j = \hat{\eta}_{u,j}(x_0; p, h) - \eta_0^{(j)}(x_0)$, $v_j = j!\hat{\beta}_j - \eta(x_0, \boldsymbol{\alpha})^\gamma (\frac{\eta_0(\cdot) - \eta(\cdot, \boldsymbol{\alpha})}{\eta(\cdot, \boldsymbol{\alpha})^\gamma})^{(j)}(x_0) - \eta(x_0, \boldsymbol{\alpha}) 1_{\{j=0\}}$ and

$$\omega_{i,j} = \binom{j}{i} (1/\eta^\gamma)^{(j-i)}(x_0, \boldsymbol{\alpha}) \eta(x_0, \boldsymbol{\alpha})^\gamma,$$

where $1_{\{j=0\}} = 1$ when $j = 0$ and $0$ otherwise.

Let $\mathbf{L}$ be a $(p+1) \times (p+1)$ matrix. For $0 \le i, j \le p$, its $(i+1, j+1)$element, denoted by $L_{i,j}$, is defined as follows. Set $L_{i,j} = 1$ when $i = j$, $L_{i,j} = 0$ when $i < j$ and

$$L_{i,j} = -\omega_{j,i} + \sum_{l=1}^{i-j-1} (-1)^{l+1} \sum_{j < k_1 < k_2 < \cdots < k_l < i} \omega_{j,k_1} \omega_{k_1,k_2} \cdots \omega_{k_l,i}$$

when $i > j$. Then $w_j = v_j - \sum_{i=0}^{j-1} \omega_{i,j} w_i = \sum_{i=0}^j L_{j,j-i} v_{j-i} = \sum_{i=0}^j L_{j,i} v_i$.

With the above notation, we have

(A.4)
$$\begin{bmatrix} \hat{\eta}_{u,0}(x_0;p,h) - \eta_0^{(0)}(x_0) \\ \hat{\eta}_{u,1}(x_0;p,h) - \eta_0^{(1)}(x_0) \\ \vdots \\ \hat{\eta}_{u,p}(x_0;p,h) - \eta_0^{(p)}(x_0) \end{bmatrix}$$

$$= \mathbf{L} \times \begin{bmatrix} 0!\hat{\beta}_0 - \eta_0(x_0) \\ 1!\hat{\beta}_1 - \eta(x_0,\boldsymbol{\alpha})^\gamma \left( \dfrac{\eta_0 - \eta(\cdot,\boldsymbol{\alpha})}{\eta(\cdot,\boldsymbol{\alpha})^\gamma} \right)^{(1)}(x_0) \\ \vdots \\ p!\hat{\beta}_p - \eta(x_0,\boldsymbol{\alpha})^\gamma \left( \dfrac{\eta_0 - \eta(\cdot,\boldsymbol{\alpha})}{\eta(\cdot,\boldsymbol{\alpha})^\gamma} \right)^{(p)}(x_0) \end{bmatrix}.$$

The above equation allows us to study the asymptotic bias and variance of $\hat{\eta}_{u,j}(x_0;p,h) - \eta_0^{(j)}(x)$ using those of $\hat{\beta}_0 - \eta_0(x_0)$ and $i!\hat{\beta}_i - \eta(x_0,\boldsymbol{\alpha})^\gamma$ $\left(\frac{\eta_0-\eta(\cdot,\boldsymbol{\alpha})}{\eta(\cdot,\boldsymbol{\alpha})^\gamma}\right)^{(i)}(x_0)$ where $1 \leq i \leq p$.

PROPOSITION 1. *Let $p - j > 0$ and suppose that conditions* (A1)–(A5) *stated in the* Appendix *are satisfied. Assume that $h = h_n \to 0$, $nh^{2p+1} \to \infty$, and $nh^{2p+3} < \infty$ as $n \to \infty$. If $x_0$ is a fixed point in the interior of* $\mathrm{supp}(f_X)$ *satisfying $\eta(x,\boldsymbol{\alpha}) \neq 0$, then*

(A.5)
$$\frac{\sqrt{nh^{2j+1}}}{\sigma_{j,j,p}(x_0;K)} \left( j!\hat{\beta}^{(j)} - \eta(x_0,\boldsymbol{\alpha})^\gamma \left( \frac{\eta_0 - \eta(\cdot,\boldsymbol{\alpha})}{\eta(\cdot,\boldsymbol{\alpha})^\gamma} \right)^{(j)}(x_0) \right.$$
$$\left. - \eta(x_0,\boldsymbol{\alpha})1_{\{j=0\}} - Bias(j) \right) \xrightarrow{\mathcal{D}} N(0,1),$$

*where the bias term $Bias(j)$ is given by $Bias_o(j)$ when $p - j$ is odd and $Bias_e(j)$ when $p - j$ is even with definitions*

$$Bias_o(j) = \frac{h^{p-j+1}}{(p+1)!} \left( \frac{\eta_0 - \eta(\cdot,\boldsymbol{\alpha})}{\eta(\cdot,\boldsymbol{\alpha})^\gamma} \right)^{(p+1)}(x_0) \eta(x_0,\boldsymbol{\alpha})^\gamma \left( \int z^{p+1} K_{j,p}(z)\, dz \right)$$
$$\times \{1 + O(h)\}$$

*and*

$$Bias_e(j) = \left\{ \int z^{p+2} K_{j,p}(z)\, dz \frac{\eta(x_0,\boldsymbol{\alpha})^\gamma}{(p+2)!} \left( \frac{\eta_0 - \eta(\cdot,\boldsymbol{\alpha})}{\eta(\cdot,\boldsymbol{\alpha})^\gamma} \right)^{(p+2)}(x_0) \right.$$
$$+ \left( \int z^{p+2} K_{j,p}(z)\, dz - j \int z^{p+1} K_{j-1,p}(z)\, dz \right) \frac{1}{(p+1)!}$$
$$\left. \times \left( \frac{\eta_0 - \eta(\cdot,\boldsymbol{\alpha})}{\eta(\cdot,\boldsymbol{\alpha})^\gamma} \right)^{(p+1)}(x_0) \frac{(\rho\eta^{2\gamma}(\cdot,\boldsymbol{\alpha})f_X)'(x_0)}{(\rho\eta^\gamma(\cdot,\boldsymbol{\alpha})f_X)(x_0)} \right\} h^{p-j+2}$$

$$\times \{1 + O(h)\}.$$

*Based on (A.4), we get the asymptotic distribution result for our estimates $\hat{\eta}_{u,j}(x_0; p, h, \boldsymbol{\alpha})$ as follows:*

$$(A.6) \qquad \frac{\sqrt{nh^{2j+1}}}{\sigma_{j,j,p}(x_0;K)}(\hat{\eta}_{u,j}(x_0;p,h) - \eta_0^{(j)}(x_0) - Bias(j)) \xrightarrow{\mathcal{D}} N(0,1).$$

*If $x_0 = x_n$ is of the form $x_0 = x_\delta + ch$ satisfying $\eta(x_0, \boldsymbol{\alpha}) \neq 0$ where $x_\delta$ is a point on the boundary of $\operatorname{supp}(f_X)$ and $c \in [-1,1]$, then (A.6) holds with $\sigma_{r,s,p}^2(x_0;K)$ and $\int z^{p+1} K_{r,p}(z)\,dz$ replaced by $\sigma_{r,s,p}^2(x_0;K,\mathcal{D}_{x_0,h})$ and $\int_{\mathcal{D}_{x_0,h}} z^{p+1} K_{r,p}(z;\mathcal{D}_{x_0,h})\,dz$.*

PROOF OF PROPOSITION 1. The result (A.6) in Proposition 1 follows from the main theorem by reading off the marginal distributions of the components of $\hat{\boldsymbol{\beta}}^*$. To calculate the asymptotic variance, we calculate the $(r+1, s+1)$ entry of $r! s! \mathbf{N}_p(\mathcal{A})^{-1} \mathbf{T}_p(\mathcal{A}) \mathbf{N}_p(\mathcal{A})^{-1}$ as

$$r! s! \sum_{k=1}^{p+1} \sum_{l=1}^{p+1} \frac{c_{r+1,k} c_{s+1,l}}{|\mathbf{N}_p(\mathcal{A})|^2} \{\mathbf{T}_p(\mathcal{A})\}_{kl} = \int_{\mathcal{A}} K_{r,p}(z;\mathcal{A}) K_{s,p}(z;\mathcal{A})\,dz,$$

where $c_{ij}$ is the cofactor of $\{\mathbf{N}_p(\mathcal{A})\}_{ij}$. The equation comes from the following argument:

$$\sum_{k=1}^{p+1} \sum_{l=1}^{p+1} c_{r+1,k} c_{s+1,l} \{\mathbf{T}_p(\mathcal{A})\}_{kl}$$

$$= \int_{\mathcal{A}} \sum_{k=1}^{p+1} \sum_{l=1}^{p+1} c_{r+1,k} c_{s+1,l} z^{k+l-2} K(z)^2\,dz$$

$$= \int_{\mathcal{A}} \left(\sum_{k=1}^{p+1} c_{r+1,k} z^{k-1} K(z)\right) \left(\sum_{k=1}^{p+1} c_{s+1,k} z^{k-1} K(z)\right) dz$$

$$= \int_{\mathcal{A}} |\mathbf{M}_{r,p}(z;\mathcal{A})| |\mathbf{M}_{s,p}(z;\mathcal{A})| K(z)^2\,dz.$$

The asymptotical results in (A.6) for $\hat{\eta}_{u,j}(x_0; p, h, \boldsymbol{\alpha})$ are easily proved by noting that its bias and variance are dominated by those of the single term $j! \hat{\beta}_j - \eta(x_0, \boldsymbol{\alpha})^\gamma (\frac{\eta_0 - \eta(\cdot, \boldsymbol{\alpha})}{\eta(\cdot, \boldsymbol{\alpha})^\gamma})^{(j)}(x_0) - \eta(x_0, \boldsymbol{\alpha}) 1_{\{j=0\}}$ based on (A.4) since $L$ is a lower triangle matrix. □

PROOF OF THEOREM 1. The result of Theorem 1 is the special case of Proposition 1 for $j = 0$. □

PROOF OF THEOREM 2. Note first that under conditions (B1)–(B5), we have $\|\hat{\boldsymbol{\alpha}} - \boldsymbol{\alpha}_0\| = n^{-1/2} O_p(1)$ by Theorem 3.2 of White (1982). This implies

that $\frac{1}{n}(Q_u(\boldsymbol{\beta}; h, x_0, \boldsymbol{\alpha}_0) - Q_u(\boldsymbol{\beta}; h, x_0, \hat{\boldsymbol{\alpha}})) \xrightarrow{P} 0$. Note that condition (A1) implies that both $Q_u(\boldsymbol{\beta}; h, x_0, \boldsymbol{\alpha}_0)$ and $Q_u(\boldsymbol{\beta}; h, x_0, \hat{\boldsymbol{\alpha}})$ are strictly concave in $\boldsymbol{\beta}$. Consequently, $\|\hat{\boldsymbol{\beta}}(x_0, \boldsymbol{\alpha}_0) - \hat{\boldsymbol{\beta}}(x_0, \hat{\boldsymbol{\alpha}})\| \xrightarrow{P} 0$ as we have obtained the asymptotic result of $\hat{\boldsymbol{\beta}}(x_0, \boldsymbol{\alpha}_0)$.

Note that $\frac{1}{n}(Q_u(\boldsymbol{\beta}; h, x_0, \boldsymbol{\alpha}_0) - Q_u(\boldsymbol{\beta}; h, x_0, \hat{\boldsymbol{\alpha}})) = O_p(1/\sqrt{n})$, $\frac{1}{n}\|\frac{\partial}{\partial \boldsymbol{\beta}} Q_u(\boldsymbol{\beta}; h, x_0, \boldsymbol{\alpha}_0) - \frac{\partial}{\partial \boldsymbol{\beta}} Q_u(\boldsymbol{\beta}; h, x_0, \hat{\boldsymbol{\alpha}})\| = O_p(1/\sqrt{n})$, $\frac{1}{n}\|\frac{\partial^2}{\partial \boldsymbol{\beta} \partial \boldsymbol{\beta}^T} Q_u(\boldsymbol{\beta}; h, x_0, \boldsymbol{\alpha}_0) - \frac{\partial^2}{\partial \boldsymbol{\beta} \partial \boldsymbol{\beta}^T} Q_u(\boldsymbol{\beta}; h, x_0, \hat{\boldsymbol{\alpha}})\|_F = O_p(1/\sqrt{n})$ for every $\boldsymbol{\beta}$ where $\|\cdot\|_F$ denotes the matrix's Frobenius norm defined as the square root of the sum of squares of each element. With consistency established above, we can consider a local compact set. By the standard argument of the Taylor expansion used for proving asymptotic normality, we get $\|\hat{\boldsymbol{\beta}}(x_0, \boldsymbol{\alpha}_0) - \hat{\boldsymbol{\beta}}(x_0, \hat{\boldsymbol{\alpha}})\| = n^{-1/2}O_p(1)$ which is faster than the convergence rate in our Theorem 1. Hence using estimated $\hat{\boldsymbol{\alpha}}$ does not affect our asymptotic convergence rates as desired. □

**Acknowledgments.** Our thanks go to an associate editor and the referees for their thorough and constructive comments.

JIANQING FAN
YANG FENG
DEPARTMENT OF OPERATIONS RESEARCH
    AND FINANCIAL ENGINEERING
PRINCETON UNIVERSITY
PRINCETON, NEW JERSEY 08544
USA
E-MAIL: jqfan@princeton.edu
        yangfeng@princeton.edu

YICHAO WU
DEPARTMENT OF STATISTICS
NORTH CAROLINA STATE UNIVERSITY
RALEIGH, NORTH CAROLINA 27695
USA
E-MAIL: wu@stat.ncsu.edu